\documentclass[10pt,a4paper,twoside,notitlepage]{article}

\usepackage{amsfonts}
\usepackage{amsmath}
\usepackage{amssymb}

\usepackage{graphicx}
\usepackage{psfrag}
\usepackage{color}

\usepackage{mathrsfs}

\newtheorem{thm}{Theorem}
\newtheorem{lem}{Lemma}[section]

\newtheorem{rem}{Remark}[section]

\newcommand{\R}{\mathbb{R}}

\newcommand{\ve}{\varepsilon}
\newcommand{\n}{\noindent}
\newcommand{\vp}{\varphi}

\newcommand{\Div}{\mathrm{div}}

\newcommand{\con}{\mathbf{C}}
\newcommand{\Lsp}{\mathbf{L}}

\newcommand{\mes}{\mathrm{meas}}

\begin{document}



\begin{center} 
\Large{{\bf Micromotions and controllability of a swimming model in an incompressible fluid governed by 2-$D$  or  3-$D$ Navier--Stokes equations}\footnote{This work  was supported in part  by  the Grant 317297 from Simons Foundation.}}\\
\large{Piermarco Cannarsa, Department of Mathematics, University of Rome ``Tor Vergata'',   Italy   and \\
 Alexander Khapalov\footnote{Corresponding author: email: khapala@math.wsu.edu},
Department of Mathematics and Statistics, \\ 
Washington State University,  USA }
\end{center}

\begin{abstract}
We study the local controllability properties of   2-$D$ and 3-$D$  bio-mimetic  swimmers  employing the change of their geometric shape to propel themselves  in  an  incompressible fluid described by Navier-Stokes equations.   It is  assumed that swimmers' bodies   consist of 
finitely many parts, identified with the fluid they occupy,  that are subsequently linked by
 the rotational  and elastic  internal forces. These  forces are explicitly described and serve as the means to affect the geometric configuration of swimmers' bodies.  Similar models were  previously investigated  in  \cite{Kh11}-\cite{Kh6}.
\end{abstract}

\section{Problem formulation and main results.}
The main goal of this paper is to study the local controllability properties of  a bio-mimetic  swimmer (see  Figures 1-8  for illustration) which makes use of its internal forces to propel itself within  a  2-$D$  or 3-$D$ incompressible fluid governed by  Navier-Stokes equations. More precisely, following   \cite{KhBook}-\cite{Kh6}, we describe swimmer's locomotion  in a fluid by the following  hybrid  nonlinear  system of two sets of partial   and ordinary differential  equations (pde/ode):
 \begin{equation}\label{eq:nse}
\left\{\begin{array}{ll}
u_t-\nu\Delta u+(u\cdot\nabla)u+\nabla p  =f & \mbox{ in } Q_T = (0,T)\times\Omega,\\
\Div \,u = 0& \mbox{ in }  Q_T,\\
u=0 & \mbox{ in } \Sigma_T =(0,T)\times\partial\Omega,\\
u(0,\cdot)=u_0 & \mbox{ in } \Omega,
\end{array}\right.
\end{equation}

\begin{equation}\label{eq:ode}
{dz_i\over dt}\,=\,{1\over \mes(S (0))}\,\int_{S(z_i(t))} u(t,x)\,dx,
\qquad 
z_i(0) = z_{i,0},
\qquad
i=1,\ldots,n, \qquad t \in (0, T).
\end{equation}
System (\ref{eq:nse})   describes the evolution of an incompressible fluid due to  Navier-Stokes equations  under the influence of the forcing term $ f (t,x)$ representing the actions of swimmer. Here,  $ x = (x_1, \ldots, x_d), d = 2, 3 $, $ \Omega$ is a bounded domain in $ \R^d$ with locally Lipschitz boundary $\partial \Omega $, $u (t,x) $ and $ p (t,x) $ are respectively  the velocity of the fluid and its pressure at point $ x$ at time $t$, and $\nu$ is the kinematic viscosity constant.  In turn, system (\ref{eq:ode}) describes the motion of  the swimmer  within $\Omega$, whose flexible body consists of $ n$   subsequently connected ``small'' sets $ S(z_i (t))$. These sets  are identified with the fluid within the space they occupy at time $t$ and are linked between themselves by the  rotational and  elastic forces as illustrated on Figures 1-7. The points $ z_i (t)$'s represent the centers of mass of the respective parts of swimmer's body. The instantaneous velocity of each part  is calculated as  the  average fluid velocity within it at time~$t$. 

\n
Below, for simplicity of notations, we will denote the sets $S(z_i(t)), i=1,\ldots,n$ also  as $S(z_i)$ or 
$S_i(t)$ and will assume the following two conditions on swimmer's body: 

\begin{description}
\item{\bf (H1)} {\em All sets $ S(z_i ), i=1,\ldots,n$ are obtained by  shifting 
the same  set $S(0)\subset\Omega$, i.e.,
$$
S(z_i)=z_i+S(0),\qquad\qquad i=1,\ldots,n, 
$$
where $S(0)$ is open and lies in a ball $B_r(0)$ of radius  $r>0$, and its center of mass is  the origin (see Figure 1). }
\end{description}

\begin{rem}\label{orient}
The main results of this paper will hold without any extra cost if we will assume that all parts of swimmer's body  are identical sets $ S(0) $ but each has its own orientation  $ S_i(0) $ in space as shown on Figure 2. In this case one  can simply change the above notations to $ S_i(z_i)= z_i+S_i(0), \;  i=1,\ldots,n$ in all respective expressions in this paper.
One can also  choose these sets to be of  distinct  shapes and sizes, in which case, however,  the respective normalizing coefficients should be added to the forcing terms to ensure that all swimmer's forces  are to be its internal forces.\end{rem}

\begin{description}\item{\bf (H2)} {\em There exist positive constants $h_0$ and ${\cal K}_S$ such that for any vector $h\in B_{h_0}(0)\setminus\{0\}$ we can find a vector $\eta=\eta\big({h}\big)$, $\;  \mid \eta \mid   = 1$ which satisfies 
\begin{equation}\label{eq:H2}
\mes(S_\Delta)_\eta^y=\int_{(S_\Delta)_\eta^y}dt\leq {\cal K}_S\,|h|
\qquad\qquad\forall~y\in\Omega_\eta\,,
\end{equation}
where $ (S_\Delta)_\eta^y$ is  the set obtained by any non-empty intersection of the set  $S_\Delta:=(h+S(0)) \,\Delta\,  S(0)=\big((h+S(0))\cup S(0)\big)\setminus \big((h+S(0))\cap S(0)\big)$ by the 
line $  \left\{y+t\eta \in \R^d  |   t \in \R,  y \in \Omega_\eta  \right\} \\ = {\mathit L}_{\eta}^y$  and $\Omega_\eta$ is the hyperplane orthogonal to $ \eta$.}
\end{description}
\bigskip
\n

\smallskip
\n
Assumption {\bf (H2)} means that  the ``thickness''  of  the set $S_\Delta$ along any line  ${\mathit L}_{\eta}^y, y \in \Omega_\eta$ parallel  to vector  $\eta$ depends uniformly Lipschitz  continuously relative to the magnitude of the shift $ h$ of the set $S(0)$  in the direction of $h$.  In the case when $ \eta (h)$ is parallel to $ h$,  {\bf (H2)} holds, e.g., for  discs and  rectangles in 2-$D$ and for balls and parallelepipeds in 3-$D$.   

 \begin{figure}[h]\label{fig:swim2d}
\centering
\includegraphics[scale=0.6]{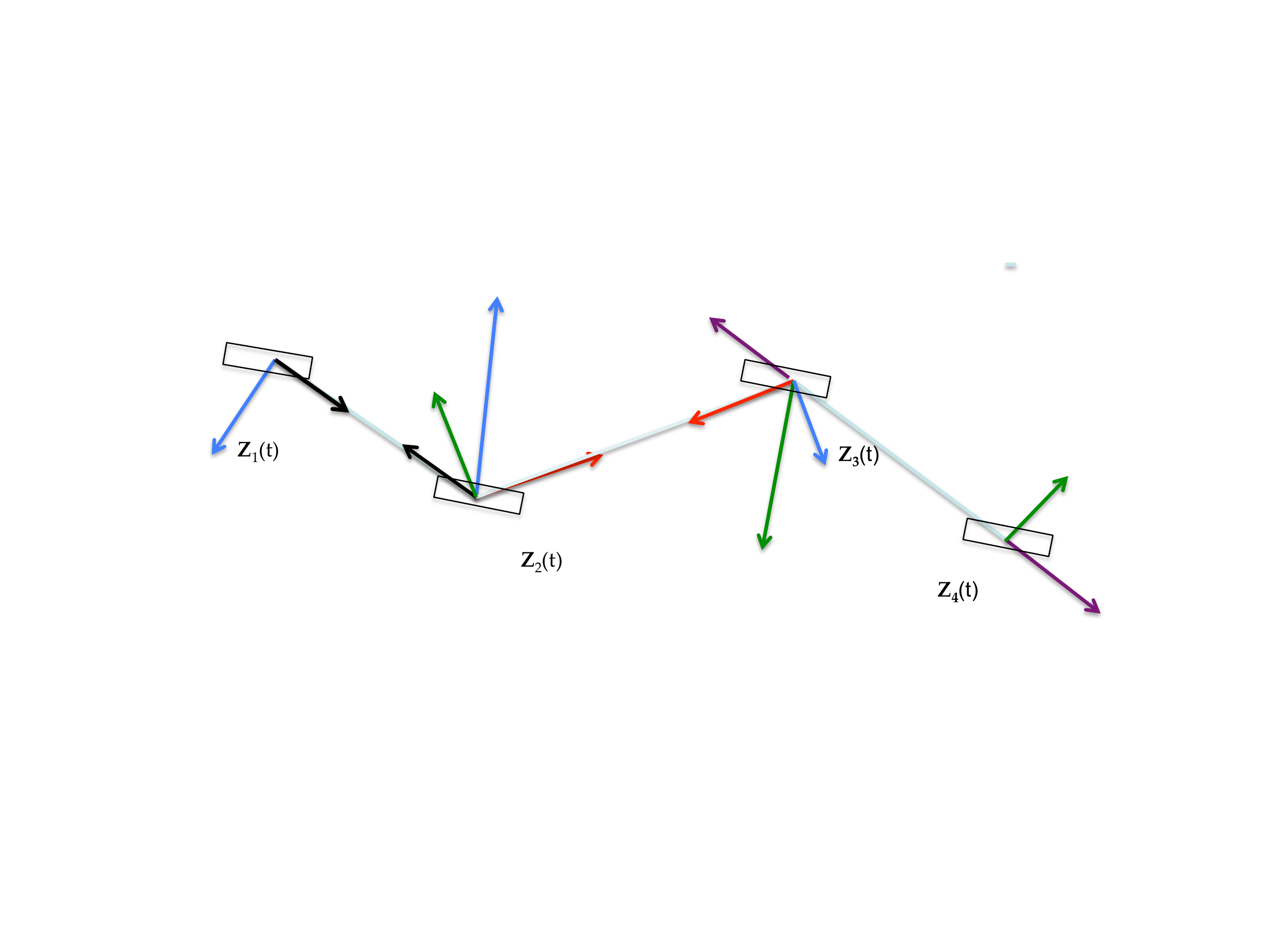} 
\caption{2-D swimmer consisting of 4 identical  rectangles of the same spatial orientation. All possible internal rotational and elastic internal forces are shown (i.e., when swimmer is not in  fluid).}
\end{figure} 

\begin{figure}[h]\label{fig:swim1d}
\centering
\includegraphics[scale=0.6]{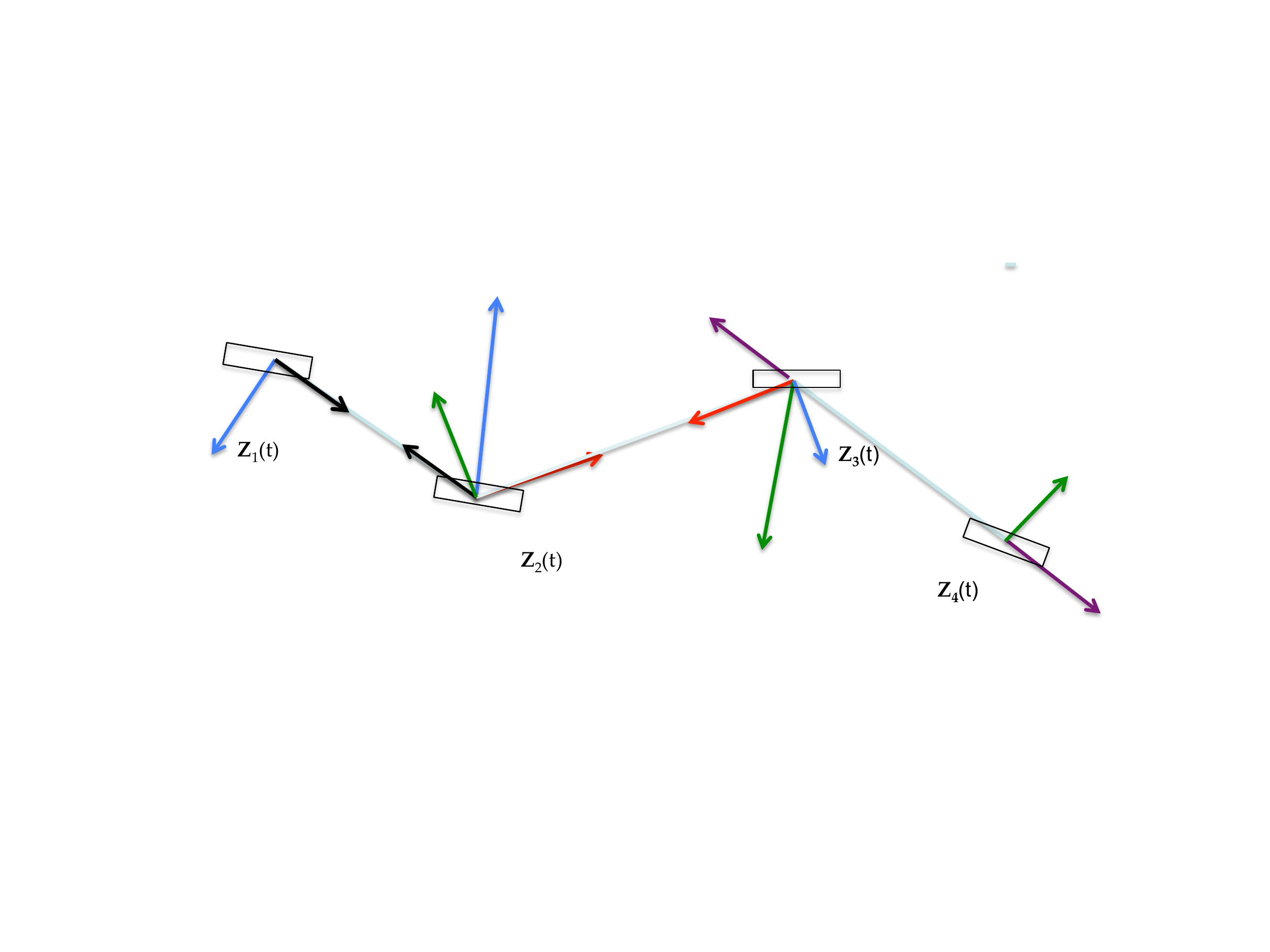} 
\caption{2-D swimmer consisting of 4 identical  rectangles which have different spatial orientation. All possible rotational and elastic internal forces are shown (i.e., when swimmer is not in  fluid).}
\end{figure}

\bigskip

\n
We assume that the forcing term $f$ in (\ref{eq:nse})  represents  the sum of \emph{rotational} and \emph{elastic}(or we can also call them {\em ``structural'') forces} generated by the swimmer (see cf. Section 12.1 in~\cite{KhBook}, \cite{Kh6}):
\begin{equation}\label{eq:forces}
f(t,x):= f_{rot}(t,x)+f_{el}(t,x).
\end{equation}
More precisely,  we assume that any of the intermediate points $ z_i, i = 2, \ldots, n-1$ can force the pair of adjacent points $ z_{i-1} $ and $ z_{i+1} $ to rotate about it, while creating,  due to Newton's 3rd Law, a counterforce acting upon $ z_i$ itself (see Figure 1):
\begin{equation}\label{rot}
f_{rot}(t,x): = \sum_{i=2}^{n-1} v_{i-1} f_{i-1} (t,x), \;\; z =(z_1, \ldots, z_n),
\end{equation}
$$
f_{i-1}(t,x) = \bigg[\xi_{i-1}(t, x)\, P_i [t] \,\big(z_{i-1}(t)-z_i(t)\big) -\xi_{i+1}(t, x)\,
{|z_{i-1}(t)-z_i(t)|^2\over |z_{i+1}(t)-z_i(t)|^2}\,
Q_i [t] \,\big(z_{i+1}(t)-z_i(t)\big)\bigg]
$$
\begin{equation}\label{rot2}
+  \xi_{i}(t, x)\,\bigg[P_i [t] \,\big(z_{i}(t)-z_{i-1}(t)\big) -\,{|z_{i-1}(t)-z_i(t)|^2\over |z_{i+1}(t)-z_i(t)|^2}\,
Q_i [t] \,\big(z_{i}(t)-z_{i+1}(t)\big)\bigg], \;\; i = 2, \ldots, n-1.
\end{equation}
In the 2-$D$ case we set  $P_i[t] = Q_i[t]  = A=\left(\begin{array}{cc}0 & 1\\ -1 & 0\end{array}\right)$. The functions $v_1,\ldots, v_{n-2} \in \Lsp^\infty(0,T)$ are {\em multiplicative  controls} (i.e., selectable parameters to control the swimming process).  

\bigskip
\n
In the 3-$D$ case, to satisfy the 3rd Newton's law, we need to make sure that the respective rotational forces acting on $z_{i-1}(t)$ and $z_{i+1}(t)$ lie in the same plane spanned by the vectors $z_{i-1}(t) - z_i(t)$ and $z_{i+1}(t) -z_i (t)$. 
 In order to achieve the continuity of these forces in time, in this paper we choose to  reduce their magnitudes to zero, when the triplet  $\{z_{i-1}(t), z_i(t), z_{i+1}(t)\}$ approaches the aligned configuration (for other options  see \cite{Kh5}). Indeed, such  configuration admits infinitely many planes containing this triplet, which makes it  an intrinsic point of discontinuity for the procedure of the choice of the rotational  plane by means of the rotational forces whose magnitudes are strictly separated from zero. Respectively, we set
  (see \cite{Kh5}-\cite{Kh6} for more details):
$$
{\bf x}\mapsto P_i[t]{\bf x}:=\big[(z_{i-1}(t)-z_i(t)) \times (z_{i+1}(t)-z_i(t))\big]\times {\bf x}\,,
$$
$$
{\bf x}\mapsto Q_i[t]{\bf x}:={\bf x}\times \big[(z_{i-1}(t)-z_i(t)) \times (z_{i+1}(t)-z_i(t))\big]\,.
$$
Note that $P_i[t]  {\bf x} =-Q_i[t] {\bf x}$ and  $|P_i[t]{\bf x}|=|Q_i[t]{\bf x}|\to 0$ for any  ${\bf x}$ when  points $z_{i-1}(t),z_i(t),z_{i+1}(t)$ converge to the  aligned configuration.

\n
In turn,  
\begin{equation}\label{elf2}
f_{el}(t,x): = \sum_{i=n}^{2n-2} v_{i-1} f_{i-1} (t,x), \;\; z =(z_1, \ldots, z_n),
\end{equation}
$$
\\
$$
\begin{equation}\label{elf}
f_{i-1}(t,x):=\xi_{i-1}(t, x)\,
\big(z_i(t) - z_{i-1}(t)\big) + \xi_{i}(t, x)\,
\big(z_{i-1}(t) - z_{i}(t)\big) ,
\end{equation}
where the functions  $v_{n-1},\ldots,v_{2n-3} \in \Lsp^\infty(0,T)$ control the distances respectively between $z_i$ and $z_{i-1}, i = 1, \ldots, n$. We set $ v = (v_1,\ldots, v_{2n-3})$.
\medskip

\n
Below,  we use the following classical notations:
\begin{itemize}
\item $\con_c^\infty(\Omega)$ denotes  the space of infinitely many times differentiable functions with compact support  in $\Omega$;
\item $H^1(\Omega)=\{\varphi | \varphi, \varphi_{x_i} \in \Lsp^2 (\Omega), i = 1, \ldots, d\} $ and $H^2(\Omega)=\{\varphi | \varphi, \varphi_{x_i}, \varphi_{x_i x_j} \in \Lsp^2 (\Omega), i,j = 1, \ldots, d\}$;
\item $H^1_0(\Omega)$ denotes the subspace of   $H^1(\Omega)$ consisting of functions vanishing on $ \partial \Omega$.
\end{itemize} 
As in \cite{Temam}, page 5, we also introduce the following  classical  vector function spaces:
$$
{\cal V}:=\{\vp\in[\con_c^\infty(\Omega]^d~|~\Div\,\vp=0\}\,, d= 2,3,
$$
$$
H:=\mathrm{cl}_{\Lsp^2}({\cal V})\,,
\qquad\qquad
V:=\mathrm{cl}_{H^1_0}({\cal V})\, =\{\vp\in [H^1_0(\Omega)]^2~|~\Div\,\vp = 0\}\,,
$$
where the symbol $\mathrm{cl}_{\Lsp^2}$ stands for  the closure with respect to the $[\Lsp^2(\Omega)]^d$-norm, and $\mathrm{cl}_{H^1_0}$ -- with respect to the $[H^1_0(\Omega)]^d$-norm. The latter is  induced  by the scalar product
$$
[ \vp,\psi ] ~:= \; \sum_{i,j=1}^d \mathop{\int}_\Omega \vp_{i x_j} \psi_{i x_j} dx.
\; \; \vp = (\vp_1,  \ldots, \vp_d), \; \psi = (\vp_1, \ldots, \vp_d).
$$

\bigskip
\n
In \cite{Kh6} we proved the following well-posedness results.

\begin{thm}[Well-posedness of  model (\ref{eq:nse})-(\ref{elf}) in the 2-$D$ case]\label{thm:wp}
Let ,  $z_{1,0},\ldots,z_{n,0} \in \Omega \subset \R^2$, $u_0 \in H$,   
and $  v_1,\ldots,v_{2n-3} \in \Lsp^\infty(0, \hat{T})$ for some $  \hat{T}>0$. Assume that 
 \begin{equation}\label{b0}
\overline{S}(z_{i,0}) \subset \Omega, \qquad\qquad
|z_{i,0}-z_{j,0}| > 2r, \qquad\qquad i, j = 1, \ldots, n,~~ i \neq j.
\end{equation} 
(Assumption (\ref{b0}) ensures that no parts of swimmer's body overlap with each other and all lie within $\Omega$.)
Then, there exists $T^*\in(0, \hat{T}]$, depending on  $u_0, z(0) = (z_{1,0},\ldots,z_{n,0})$ and the $\Lsp^\infty(0,  \hat{T})$-norms of $ v_j$'s, such that system~(\ref{eq:nse})-(\ref{b0}) admits a unique solution 
$$
(u, z) \in 
\bigg(C([0,T^*]; H) \bigcap \Lsp^2 (0, T^*; V) \bigg)  \times [C([0,T^*]\,;\R^2)]^n, 
$$
and 
 \begin{equation}\label{body2}
\overline{S}(z_i(t)) \subset \Omega, \qquad
|z_i(t)-z_j(t)| > 2r, \qquad i, j = 1, \ldots, n,~~ i \neq j \, \qquad \forall  t \in [0, T^*].
\end{equation}
\end{thm}
The equation (\ref{eq:nse}) in the above is understood  in the sense of the following identity:
$$
\int_\Omega u (x, t) \phi(t) dx -  \int_\Omega u (x, 0) \phi(t) dx - \int_0^{t} \int_\Omega u \cdot \phi_t dx dt  \; =  
$$
\begin{equation}\label{identityu}
=  -  \int_0^{t}  \nu\,[ u (\tau, \cdot) ,\phi (\tau, \cdot) ] dt \; + \;  \int_0^{t} \int_\Omega \bigg(- (u \cdot \nabla)  u +\, P_H f_j(\tau,x;h)  \bigg) \phi dx d\tau, \;\; t \in [0, T^*],
\end{equation}
where   $ P_H: [\Lsp^2(\Omega)]^d \rightarrow H $  denotes  the projection operator from $[\Lsp^2(\Omega)]^d$ onto $H$ and $ \phi $ is any function such that $ \phi \in  \Phi = \{\phi \in  \Lsp^2 (0, T^*; V) $, $ \phi_t \in  \Lsp^2 (0,T^*; H)\}$. 
The complementing $ \nabla p  = f -u_t + \nu \Delta u - (u \cdot \nabla )u $ is understood in the sense of distributions.

\begin{thm}[Well-posedness of  model (\ref{eq:nse})-(\ref{elf}) in 2-$D$ and 3-$D$ cases]\label{thm:wp23d}
Let $\partial \Omega$ be of class $ C^2$,    $z_{1,0},\ldots,z_{n,0} \in \Omega \subset \R^d, \; d = 2, 3$, $u_0 \in  V$,   
and $  v_1,\ldots,v_{2n-3} \in \Lsp^\infty(0, \hat{T})$ for some $  \hat{T}>0$. Assume that 
(\ref{b0}) holds.
Then, there exists $T^*\in(0, \hat{T}]$, depending on  $u_0, z(0) = (z_{1,0},\ldots,z_{n,0})$ and the $\Lsp^\infty(0,  \hat{T})$-norms of $ v_j$'s, such that system~(\ref{eq:nse})-(\ref{b0}) admits a unique solution 
$$
(u, z) \in  
C([0,T^*]; V)  \times [C([0,T^*]\,;\R^d)]^n, \;\;\;\; \nabla p \in [\Lsp^{2} (Q_{T^*})]^d
$$
($ \nabla p (t, \cdot) \in H^\perp$,  the orthogonal complement of $ H$ in $ [\Lsp^2 (\Omega)]^d$) and 
 \begin{equation}\label{body2}
\overline{S}(z_i(t)) \subset \Omega, \qquad
|z_i(t)-z_j(t)| > 2r, \qquad i, j = 1, \ldots, n,~~ i \neq j \, \qquad \forall  t \in [0, T^*].
\end{equation}
\end{thm}

\begin{rem}\label{fp}
The  argument of  \cite{Kh6} makes use of  Schauder's fixed point theorem. In particular, we showed that the sequence of uncoupled mappings corresponding to the uncoupled version of system 
(\ref{eq:nse})-(\ref{b0}), namely:  $ \; w $ (in place of $ u$) $ \rightarrow z \rightarrow f \rightarrow u$ are continuos with respect to the norms $ \Lsp^2 (0,T^*; V) \rightarrow  [C([0,T^*]\,;\R^2)]^n \rightarrow \Lsp^2 (0,T^*; [\Lsp^2 (\Omega)]^d) \rightarrow \Lsp^2 (0,T^*; V)$ and their product is a compact operator with the unique fixed point $ u$.
\end{rem}

\begin{rem}[Some useful estimates]\label{nseest} Solution to (\ref{eq:nse}) satisfies the following estimates (see \cite{Lad70},  Lemma 9, p. 194, (55); \cite{Lad69},  \cite{Temam} \cite{Kh6}).
\begin{itemize}
\item
Under the assumptions of Theorems  \ref{thm:wp} and \ref{thm:wp23d}:
\begin{equation}\label{th1}
\|u\|_{C (0,T^*\,; H)} + \|u\|_{\Lsp^2(0,T^*\,;V)}\leq   \; L_{\bf S_1} \left(\,\|u_0\|_{H}+ \,\|f\|_{\Lsp^{2}(0,T^*\,;[\Lsp^2 (\Omega)]^d)} \right)
\end{equation}
for some constant $ L_{\bf S_1}$. 
Let us also recall here  the following estimate from Theorem 11 in \cite{Lad69}, see estimates (45) and (48) on pp. 170-171 (see also \cite{Kh6}) 
\begin{equation}\label{th11}
\| u_{(1)}- u_{(2)}  \|_{C([0,T^*]; H)} \; + \| u_{(1)}- u_{(2)}  \|_{\Lsp^2 (0,T^*; V)}\leq \; D^* \;  \| f_{(1 )} - f_{(2)}  \|_{\Lsp^{2}(0,T^*\,;[\Lsp^2 (\Omega)]^d)} ,
\end{equation}
where $ u_{(m)} $ is the solution to (\ref{eq:nse})-(\ref{elf}) for $ f = f_{(m )}, m = 1, 2$ and $ D^*(s)$ is a nondecreasing function of  $ \max_{m = 1, 2}  \{ \| u_{(m)}  \|_{\Lsp^2 (0,T^*; V)})\}$.
\item
Under the assumptions of Theorems  \ref{thm:wp} and \ref{thm:wp23d}:
\begin{equation}\label{fje}
\| v_j f_j \|_{[\Lsp^\infty (Q_{T^*})]^d} \; \leq \; C_{\Omega, r} | v_j|, \;\; d = 2, 3,
\end{equation}
for some constant $ C_{\Omega, r} > $ depending on $\Omega$ and $ r$.
\item
In tern, under the assumptions of Theorem \ref{thm:wp23d}:
\begin{equation}\label{th12}
\|u\|_{C (0,T^*\,; V)} \;  \leq   \; L_{\bf S_2} \left(\,\|u_0\|_{V}+ \,\|f\|_{\Lsp^{2}(0,T^*\,;[\Lsp^2 (\Omega)]^d)} \right)
\end{equation}
for some constant $ L_{\bf S_2}$. 
\end{itemize}
\end{rem}

\n
\begin{description}\item{\bf (H3)} 
{\em Everywhere below we assume that the initial datum  $ u_0   $ is fixed.}
\end{description}

\bigskip
\n 
Our first main  result  describes  the micromotions of the swimmer in 2-$D$ and 3-$D$ models (\ref{eq:nse})-(\ref{b0}) in terms of projections  of its internal forces at the initial moment on $H$. 

\begin{thm}[Swimmer's micromotions] \label{mm} Under the assumptions of Theorems \ref{thm:wp} and  \ref{thm:wp23d}, if we set 
$
v_j = h a_j \in \R, \; \sum_{j = 1}^{2n - 3} a_j^2 = 1, t \in (0, T^*), \; \| (v _1, \ldots, v_{2n-3})\|_{\R^{2n-3}} = | h | \leq 1$ , then 
\begin{equation}\label{zijdirall}
z_i (t) \;  =  z_i(0) \; +    \;{h t^2 \over {2 \mes(S(0))}} \, \sum_{j = 1}^{2n - 3}  a_j \int_{S(z_i(\tau;0))} (P_H f_j^T (0, \cdot))(x) dx + ht^2 \rho (t)  + h \zeta (h,t), \;\; t \in  [0, T^*], 
\end{equation}
where $ \| \rho \|_{[C[0,t]]^d}  =   O(t) $,  $  \| \zeta (h, \cdot) \|_{[C[0,t]]^d}  = \;  O(h)$ ($d = 2, 3$)  and are defined by $ u_0$ and $ f_j(0, \cdot)$. (Here and below $ O(p)$ denotes a real-valued function that tends to 0 as $ \R \ni p \rightarrow 0$.)
\end{thm}
The main controllability results of  this paper are as follows.

\begin{thm}[Local swimming controllability of $z_i$'s: 2-$D$ case] \label{thm:controllability} 
Given $ u_0  \in  H$, 
under the assumptions of Theorem \ref{thm:wp},  let $ (u^*, z^* =  (z_1^*,\ldots,z^*_n)) $ be the solution to (\ref{eq:nse})-(\ref{elf}) generated by the  zero controls $v_1=\ldots=v_{2n-3} = 0$  on some interval $[0, T^*]$. (Note: due to Remark 1.2 and (\ref{body2}),  the curves $ z^*_i (t), t \in [0, T^*], i = 1, \ldots, n $ lie in $ \Omega$ along with some their neighborhoods: 
\begin{equation}\label{conv}
B_\mu (z_i^*(t)) \subset \Omega, \;\; t \in [0, T^*], \;\; i = 1, \ldots, n,
 \end{equation}
where  $ B_\mu (z^*(t)) $ is the ball in $ \R^2$ with center at $ z_i^* (t) $ of radius $ \mu>0$.)   Let for some $ i \in \{1, \ldots, n\}$ and $k,l \in\{1,\ldots,2n-3\}$ the  vectors
\begin{equation}\label{eq:indep_vect}
 \int_{S(z_i(0))} (P_H  f_k (0, \cdot))(x)\,dx,  \; 
\int_{S(z_i(0))} (P_H  f_l (0, \cdot))(x)\,dx 
\end{equation}
 be linearly independent. Then there exist $T=T(i,k,\ell)\in\,(0,T^*]$ and $\ve=\ve(i,k,\ell)>0$ such that
$$
B_\ve(z_i^*(T))\subseteq
\Big\{
z_i(T)~|~v_k,v_\ell\in\R, \;\;  {\rm while} \;\;  v_j=0 \;\; {\rm for} \;\;  j=1,\ldots,2n-3,~j\neq k,\ell
\Big\}\,.
$$
\end{thm}
In other words, under the conditions of Theorem \ref{thm:controllability}, the point $ z_i$ can be steered on some  time-interval $[0, T]$  from its initial position $ z_{i, 0} = z_i^*(0 )$ to any point within the ball $ B_\ve(z_i^*(T))$ of radius $ \ve>0$ with center  at  the endpoint $z_i^*(T)$  of  the ``drifting'' trajectory $ z_i^* (t), t \in [0, T]$. We will also show in our proofs below that {\em this can be achieved  merely by constant controls} $ v_i$'s. 

\bigskip
\n
At no extra cost (making use of (\ref{zijdirall}) instead of (\ref{zikldir})), we will have the following result for the motion of the center of mass of our swimmer.

\begin{thm}[Local swimming locomotion: 2-$D$ case] \label{thm:ccm} 
Let in Theorem \ref{thm:controllability} condition (\ref{eq:indep_vect})  is replaced with the following: 
\begin{equation}\label{eq:indep_vectc}
\sum_{i = 1}^n  \int_{S(z_i(0))} (P_H  f_k (0, \cdot))(x)\,dx,  \; 
\sum_{i = 1}^n  \int_{S(z_i(0))} (P_H  f_l (0, \cdot))(x)\,dx 
\end{equation} 
are linearly independent. Then the result of Theorem \ref{thm:controllability} holds with respect to the swimmer's center of mass $ z_c = \frac{1}{n} \sum_{i = 1}^n z_i (t)$, namely:
$$
B_\ve(z_c^*(T))\subseteq
\Big\{
z_c(T)~|~v_k,v_\ell\in\R,\;  {\rm while} \;\;  v_j=0 \;\; {\rm for} \; j=1,\ldots,2n-3,~j\neq k,\ell\Big\}, \; z_c^* = \frac{1}{n} \sum_{i = 1}^n z_i^* (t).
$$
\end{thm}

\begin{thm}[Local controllabity in 3-$D$] \label{thm:c3d} Given $ u_0  \in  V$, 
under the assumptions of Theorem \ref{thm:wp23d}, the results of Theorem \ref{thm:controllability}  can be extended to the case of 3-$D$ swimming model  (\ref{eq:nse})-(\ref{elf}), assuming that  three  controls  $v_k, v_l $ and $ v_m$ are active (i.e., in place of  two as in Theorem \ref{thm:controllability}).
\end{thm}

\begin{thm}[Local locomotion  in 3-$D$] \label{thm:c3dl} Given $ u_0  \in  V$, 
under the assumptions of Theorem \ref{thm:wp23d}, the results of Theorem \ref{thm:ccm}   hold for the case of 3-$D$ swimming model  (\ref{eq:nse})-(\ref{elf}) for three  active controls  $v_k, v_l $ and $ v_m$  (i.e., in place of  two as in Theorem \ref{thm:ccm}).
\end{thm}

\bigskip
\n
The main idea of our proofs below is  to show that each of the mappings 
\begin{equation}\label{ifth}
\R^2 \ni (v_k, v_l) \; \rightarrow \; z_i (T) \in \R^2, \;\;\; \R^3 \ni (v_k, v_l, v_m) \; \rightarrow \; z_i (T) \in \R^3,
\end{equation}
associated with   Theorems \ref{thm:controllability} and \ref{thm:c3d}, 
considered on some (open) neighborhood of the origin, is 1-1 and its range contains an open neighborhood of $ z_i^* (T)$ for some $ T>0$.
To this end, we intend to study the invertibility properties of the respective $[2\times2]$- and $[3\times3]$-matrices:
$$
\left( \frac{d z_i (T) }{dv_k} \mid_{v_j's = 0} ,  \frac{d z_i (T) }{dv_l} \mid_{v_j's = 0} \right)
$$
\begin{equation}\label{mtx}
\left( \frac{d z_i (T) }{dv_k} \mid_{v_j's = 0} ,  \frac{d z_i (T) }{dv_l} \mid_{v_j's = 0},  \frac{d z_i (T) }{dv_m} \mid_{v_j's = 0} \right).
\end{equation}
In the above and anywhere below the subscript $v_j's = 0$ indicates that the corresponding expressions are calculated for $ \; v_j = 0, j = 1, \ldots, 2n-3$.

\bigskip
\n
The remainder of the paper is organized as follows. Sections 2 and 3 deal with detailed proofs of auxiliary results in the 2-$D$ case. Namely, in Section 2  we will describe the derivatives ${\partial u\over \partial v_j} |_{v_j's=0}, j = 1, \ldots, 2n-3$ as solutions to  some linear system of partial differential equations. Then in Section 3 we will show that  the derivatives  ${\partial z_i\over \partial v_j} |_ {v_j's = 0}, i = 1, \ldots, n, \; j = 1, \ldots, 2n-3$ satisfy a system of integral   Volterra equations. In Section 4 we will prove Theorem \ref{mm} and the main controllability results for the 2-$D$ case. In Section 5 we show how they can be extended to the 3-$D$ case. In Section 6 we discuss illustrating examples.

\bigskip
\n
{\bf Prior related results on swimming controllability}.
Local controllability results  similar to Theorems \ref{thm:controllability} and \ref{thm:ccm}  were obtained in \cite{Kh7} (see also \cite{KhBook}, Ch. 14) for the case of an incompressible fluid governed by the non-stationary Stokes equations and when the elastic forces were described by ``uncontrollable'' Hooke's Law. 
The proofs in \cite{Kh7} were based on the Inverse Function Theorem for the 2-$D$ mapping (\ref{ifth}) and employed the linearity of the fluid equations to represent the velocity of fluid in the  form of implicit Fourier series expanded along the associated  set of eigenfunctions. In this paper we  follow  the  general strategy of  \cite{Kh7}. However, the nonlinearity  of Navier-Stokes equations requires a principal  modification of  this strategy and its setup.  In particular, in (\ref{elf2})-(\ref{elf}) we consider controlled elastic forces instead of Hooke's Law as in \cite{Kh7}. For such forces in \cite{KhBook}, Ch. 15  we obtained some global controllability results for the case a swimmer applying a rowing-type motion in a fluid governed  by the non-stationary Stokes equations (see also \cite{Kh4} for the 3-$D$ case).

\begin{rem}[Swimming controllability in the framework of ODE's]\label{swcode}
A number of attempts were made to study controllability of  various `` swimmers''
in the context of swimming models in the framework of ODE's, see, e.g., 
 Koiller et al.  \cite{Koi}  (1996); McIsaac and
Ostrowski \cite{McI}  (2000); Martinez and Cortes  \cite{Mar}  (2001);
Trintafyllou et al.  \cite{Tri} (2000); 
Alouges et al. \cite{Alo} (2008), Sigalotti and  Vivalda  \cite{Sig}  (2009),
and the references therein.  
\end{rem}


\section{Derivatives ${\partial u\over \partial v_j} |_{v_j's=0} $: 2-$D$ case}\label{duv}
As suggested by the representaion of matrices in (\ref{mtx}),  we intend to  evaluate the derivatives ${d\over d v_j}\, z_i(t)$ assuming that $ v_j$'s are independent   variables (real numbers) in (\ref{eq:nse})-(\ref{eq:ode}). As the 1st step in this direction, in this section we will study  derivatives ${\partial u\over \partial v_j} |_{v_j's=0}$.

\n
Fix any $ j \in \{1, \dots, 2n-3\} $ and  assume that in  (\ref{eq:nse})-(\ref{elf})  
\begin{equation}\label{h1}
v_j = h \in \R, \;\; | h | \leq 1,\;\; v_m = 0, \;\; m \neq j,
\end{equation}
denoting respectively  in this case:
$$
 z_i (t)  =  z_i (t;h),  \;\; 
 z (t)  =  z (t;h),  \;\; f(t,x) = f(t, x;h), \;\;  f_j (t,x) = f_j (t,x;  h), \;  
    p (t,x)  = p(t,x;  h), 
    $$
\begin{equation}\label{fh1}
  u(t, x)= u (t, x; h) =  u_h (t, x), \;\; u_* (t, x) =  u (t,x; 0), \;\;  w_h=\, \frac{u_h - u_*}{h}.
\end{equation}
 We will study the behavior of $w_h$ as $ h $ tends to zero.
  Then, (\ref{eq:nse}) yields:
\begin{equation}\label{eq:wh}
\left\{\begin{array}{ll}
w_{ht} =\nu\,\Delta w_h - (u_*\cdot\nabla) w_h-(w_h \cdot\nabla) u_h  \\
+ \; \, f_{j}(\cdot,\cdot;h)- {1 \over h} \nabla (p (\cdot, \cdot; h) - p (\cdot, \cdot; 0)) & \mbox{ in } (0,T^*)\times\Omega,\\
\Div \,w_h = 0& \mbox{ in } (0,T^*)\times\Omega,\; {\rm i.e.}, \; w_h (t, \cdot) \in H, \\
w_h=0 & \mbox{ in } (0,T^*)\times\partial\Omega,\\
w_h(0,\cdot)=0 & \mbox{ in } \Omega.
\end{array}\right.
\end{equation}

\n
By Theorem \ref{thm:wp}, (\ref{rot2}), (\ref{elf}) and  
due to continuous embedding (for $\Omega \subset \R^2 $, see Remark \ref{l4} below)
\begin{equation}\label{eml4}
C ([0, T^*]; [\Lsp^2 (\Omega)]^2) \bigcap \Lsp^2([0, T^*; V)  \; \subset [\Lsp^4 (Q_{T^*})]^2= [\Lsp_{4,4} (Q_{T^*})]^2,
\end{equation}
we have:
\begin{equation}\label{u*est}
u_*, u_h, w_h  \in C ([0, T^*]; H) \bigcap \Lsp^2 (0, T^*; V)\bigcap [\Lsp^4 (Q_{T^*})]^2,   \;\;  f_j (\cdot, \cdot; h) \in[ \Lsp^\infty (Q_{T^*})]^2.
\end{equation}

\begin{rem}\label{l4}
 In the above we used  estimate  (3.4) in  \cite{Lad2}, page 75, namely:
\begin{equation}\label{4q}
\| \psi \|_{\Lsp^4(Q_{t})} \; \leq \; \beta \bigg( \max_{\tau \in [0, t]} \| \psi (\tau, \cdot) \|_{\Lsp^2(\Omega )}  + \| \nabla \psi \|_{[\Lsp^2(Q_{t})]^2} \bigg).
\end{equation}
\end{rem}

\bigskip
\n
We claim that Theorem 1.1 in \cite{Lad2}, pages 573-574  on well-posednes of general parabolic systems (see Remark \ref{ch} below)
implies that (\ref{eq:wh}) admits a unique solution of regularity described in (\ref{u*est}) 
and for some constant $C^* > 0$ the following estimate holds: 
$$
\| w_h \|_{C ([0, T^*]; [\Lsp^2 (\Omega)]^2) \bigcap \Lsp^2([0, T^*; V)} \mathop{=}^\Delta  \; \max_{t \in [0, T^*]} \| w_h (t, \cdot) \|_{[\Lsp^2 (\Omega)]^2}
$$
\begin{equation}\label{whest} 
 + \| w \|_{\Lsp^2 (0,T^*; V)} \; \leq C^*  \| P f_{j}(\cdot,\cdot;h)  \|_{[\Lsp_{2, 1}(Q_{T^*})]^2} \leq C^*  \| f_{j}(\cdot,\cdot;h)  \|_{[\Lsp_{2, 1}(Q_{T^*})]^2}. 
\end{equation}

\n
Indeed, the proof of this theorem  is based on Galerkin methods with test functions 
$$
 \phi  \in \Lsp^2 (0, T^*; [H_0^1 (\Omega)]^2) \bigcap 
\{ \phi (\cdot, x) \in [H^1 (0, T^*)]^2 \;{\rm\; a.e. \; in} \;  \Omega\}, 
$$ 
see \cite{Lad2}. However, in the case of the special mixed problem (\ref{eq:wh}), including the extra condition that $ \Div \, w_h = 0 $, we are dealing with  $ w_h (t, \cdot) $ that lie in $ V$ for almost all $t$, see (\ref{u*est}). Therefore, $ w_h$ can be represented  as a Fourier series expanded {\em only} along the eigenfunctions $ \{\omega_k\}_{k = 1}^\infty, \omega_k \in  V, k = 1, \ldots $ of the spectral problem associated with (\ref{eq:nse}), forming a complete orthogonal basis in $ V $ and orthonormal in $ H$ (\cite{Lad69}), namely, in the following form:
\begin{equation}\label{Galerkin}
w_h (t,x) = \sum_{k=1}^\infty c_k (t) \omega_k (x); \;\;\;\; \nu \Delta \omega_k = \lambda_k \omega_k + \nabla p_k,   \; \Div \, \omega_ k = 0 \;\; {\rm in} \; \Omega, \;\; k = 1, \ldots.
\end{equation}
The equation for $ \omega_k$'s  is understood in the sense of  identity (see also (\ref{identitywh}))
$$
-  \nu\,[ w_k ,\phi  ]  \; = \;  \lambda_k  \int_\Omega\omega_k  \phi dx d\tau \;\;\;\; \forall \phi \in V.
$$
In other words, (\ref{eq:wh})  is equivalent to the following identity obtained as the difference of identities (\ref{identityu}) in the cases when   $ u = u_h $ and $ u = u_*$ and then divided by $h$ (compare to \cite{Lad2}, p. 572):
$$
\int_\Omega w_h (x, t) \phi(t) dx - \int_0^{t} \int_\Omega w_h \cdot \phi_t dx dt  \; 
=  
$$
\begin{equation}\label{identitywh}
=  -  \int_0^{t}  \nu\,[ w_h (\tau, \cdot) ,\phi (\tau, \cdot) ] dt \; + \;  \int_0^{t} \int_\Omega \bigg(- (u_* \cdot \nabla) w_h-(w_h \cdot \nabla) u_h+\, P_H f_j(\tau,x;h)  \bigg) \phi dx d\tau, \;\; t \in [0, T^*], 
\end{equation}
where $ \phi $ is any function such that $ \phi \in  \Phi = \{\phi \in  \Lsp^2 (0, T^*; V) $, $ \phi_t \in  \Lsp^2 (0,T^*; H)\}$.

 \n
 The derivation of (\ref{whest}) in \cite{Lad2} is based on the classical form  of identity  (\ref{identitywh}), namely,  with any $ \phi \in  \Lsp^2 (0, T^*; [H_0^1 (\Omega)]^2),  \phi_t  \in \Lsp^2 (Q_{T^*}) $, which in this case will be applied for 
  $\phi  \in \Phi$, and, thus,  is the same as just to use the last line in (\ref{identitywh}) from the start. 
  
  \n
  
\begin{rem}\label{ch} Let us recall the following results from \cite{Lad2}.
\begin{itemize}
\item
Theorem 1.1 in \cite{Lad2}, pages 573 requires that the  squared 1-$D$ components of the 2-$D$ vector-function  $u_*$    and 1-$D$ components of the $2 \times 2$ matrix-function $ \nabla u_h$ (as the coefficients in (\ref{eq:wh})) are elements of the space $\Lsp_{q,r} (Q_{T^*})$, where 
$$
\| \psi \|_{\Lsp_{q,r} (Q_{T^*})} \; = \; \bigg(\int_0^{T^*} \bigg(\int_\Omega | \psi |^q dx\bigg)^{r/q} dt\bigg)^{1/r}, \;\; \frac{1}{r} + \frac{1}{q} = 1, \;\;  q \in (1, \infty], \; r \in [1, \infty),
$$
while the free term $f_{j}(\cdot,\cdot;h) $ lies in $\Lsp_{q_1,r_1} (Q_{T^*})$ (due to (\ref{identitywh}) we ignore the term $ - {1 \over h} \nabla (p (\cdot, \cdot; h) - p (\cdot, \cdot; 0))$ here), where $ \frac{1}{r_1} + \frac{1}{q_1} = 1 + \frac{1}{2}, \;  
q_1 \in (1, 2], \; r_1 \in [1, 2)$. We can select $ r_1 = 1, q_1 = 2$ and  $ r = 2 = q $.
\item
Constant $C^*$ can be selected to be dependent only on 
$\| u_* \|_{C ([0, T^*]; [\Lsp^2 (\Omega)]^2) \bigcap \Lsp^2([0, T^*; V)}$ or $ \| u_*\|_{[\Lsp^4 (Q_{T^*})]^2}$,  and $\| \nabla u_h \|_{[[\Lsp^2 (Q_{T^*})]^2]^2}, \, | h | \leq 1$, see \cite{Lad2}, pages 573-574 and (\ref{h1}).
\item
Condition $ \Div \, w_h = 0 $ is not required in Theorem 1.1 in \cite{Lad2}, page 573.
\item
Note that  $w_h$ also satisfy the regularity of solutions to (\ref{eq:nse}) described in Theorem \ref{thm:wp}. 
\end{itemize}
\end{rem}  
  \n
Based on the above discussion, we can refine   (\ref{whest}) as follows: 
$$
  \| w_h \|_{C ([0, T^*]; [\Lsp^2 (\Omega)]^2) \bigcap \Lsp^2([0, T^*; V)}  \; \leq C^*  \|  f_{j}(\cdot,\cdot;h)  \|_{[\Lsp_{2, 1}(Q_{T^*})]^2}
$$
\begin{equation}\label{whest2} 
 \leq C^* T^* \, {\rm meas}^{1/2} \, (\Omega) \,   \|  f_j (\cdot, \cdot;h) \|_{[\Lsp^\infty(Q_{T^*})]^2} \; 
\leq \; C^* T^* \, {\rm meas}^{1/2} \, (\Omega) \,   C_{\Omega, r} ,  \;\; j = 1, \ldots, 2n-3.
\end{equation}

\bigskip
\n
Introduce the following linear system:
\begin{equation}\label{eq:wj}
\left\{\begin{array}{ll}
\bigg({\partial u\over \partial v_j} |_{v_j's=0} \bigg)_t \; = \; \nu\,\Delta \bigg( {\partial u\over \partial v_j} |_{v_j's=0} \bigg) - (u_*\cdot\nabla) \bigg( {\partial u\over \partial v_j} |_{v_j's=0} \bigg) \\
-(\bigg( {\partial u\over \partial v_j} |_{v_j's=0} \bigg) \cdot\nabla) u_*+\, P_H f_j (\cdot, \cdot;0)  - \nabla p_j & \mbox{ in } (0,T^*)\times\Omega,\\
\Div \,\bigg( {\partial u\over \partial v_j} |_{v_j's=0} \bigg) = 0& \mbox{ in } (0,T^*)\times\Omega, \\
 & \; {\rm i.e.}, \; \bigg( {\partial u  \over \partial v_j} |_{v_j's=0} \bigg) (t, \cdot)  \in H, \\
\bigg( {\partial u\over \partial v_j} |_{v_j's=0} \bigg)=0 & \mbox{ in } (0,T^*)\times\partial\Omega,\\
\bigg( {\partial u\over \partial v_j} |_{v_j's=0} \bigg)(0,\cdot)=0 & \mbox{ in } \Omega,
\end{array}\right.
\end{equation}
where (in the sense of distributions, see also Theorem \ref{thm:wp}), 
$$
\nabla p_j  = 
$$
$$
-  \bigg({\partial u\over \partial v_j} |_{v_j's=0} \bigg)_t +  \nu\,\Delta \bigg( {\partial u\over \partial v_j} |_{v_j's=0} \bigg) - (u_*\cdot\nabla) \bigg( {\partial u\over \partial v_j} |_{v_j's=0} \bigg) 
-(\bigg( {\partial u\over \partial v_j} |_{v_j's=0} \bigg) \cdot\nabla) u_*+ P_H f_j (\cdot, \cdot;0).
$$

\bigskip

\begin{rem}[On understanding system (\ref{eq:wj})]\label{rgal}
Due to incompressibility (``divergence-free'') condition $$ \Div \,\bigg( {\partial u\over \partial v_j} |_{v_j's=0} \bigg) = 0,$$ similar to (\ref{Galerkin}), we can represent solution to (\ref{eq:wj}) as the series 
\begin{equation}\label{Galerkinwj}
\bigg({\partial u\over \partial v_j} |_{v_j's=0} \bigg)(t,x) = \sum_{k=1}^\infty d_k (t) \omega_k (x).
\end{equation}
Then the argument of the classical theory of parabolic pde's (\cite{Lad2}, Chapters VII and III) can be applied to the ``cut-off''  form (\ref{Galerkinwj}) exactly as it is applied in the case when such condition is absent.  Respectively, exactly as in the aforementioned classical theory, making use of the identity like (\ref{identitywh}), we can derive the existence of solution to (\ref{eq:wj})  in $C ([0, T^*]; [\Lsp^2 (\Omega)]^2) \bigcap \Lsp^2([0, T^*; V)$ satisfying (\ref{est:wj}).
\end{rem}

\begin{lem}\label{lem:wj}  Derivatives 
$$ 
\lim_{h\rightarrow 0} w_h \mathop{=}^\Delta  {\partial u\over \partial v_j} |_{v_j's=0}, \;\;\;\; j = 1, \ldots, 2n -3, 
$$
where the limit is taken with respect to the $C ([0, T^*]; [\Lsp^2 (\Omega)]^2) \bigcap \Lsp^2([0, T^*; V)$-norm, exist  as unique solutions to  (\ref{eq:wj})  and
\begin{equation}\label{wprop}
{\partial u\over \partial v_j} |_{v_j's=0} \; \in \; C ([0, T^*]; H) \bigcap \Lsp^2 (0, T^*; V) \bigcap [\Lsp^4 (Q_{T^*})]^2.
\end{equation} 
 As a particular case of  (\ref{whest2}), the following estimates hold:
$$
\|  {\partial u\over \partial v_j} |_{v_j's=0} \|_{C ([0, T^*]; H) \bigcap \Lsp^2 (0, T^*; V)} \; 
 $$
\begin{equation}\label{est:wj} 
\leq C^* T^* \, {\rm meas}^{1/2} \, (\Omega) \,   \|  f_j (\cdot, \cdot;0 ) \|_{[\Lsp^\infty(Q_{T^*})]^2} \; 
\leq \; C^* T^* \, {\rm meas}^{1/2} \, (\Omega) \,   C_{\Omega, r} ,  \;\; j = 1, \ldots, 2n-3. 
\end{equation}
where $C^*$ can be selected to be dependent only on $ \| u_*\|_{[\Lsp^4 (Q_{T^*})]^2}$ and  $\| \nabla u_* \|_{[\Lsp^2 (Q_{T^*})]^2}$.

\end{lem}

\n
{\bf Proof:}  

\n
{\bf Step 1.}
We can  use here an adoptation  of the argument of Theorem 4.5 in \cite{Lad2}, page 166 (on continuous dependence of solutions to parabolic pde's on  coefficients and free terms) to the case of systems of linear parabolic pde's   along Remark \ref{rgal}. Namely, denote 
$$
W_h  = w_h - {\partial u\over \partial v_j} |_{v_j's=0} .
$$
Then, we will have the following identity for $ W_h$ from (\ref{identitywh}):
$$
\int_\Omega W_h (x, t) \phi(t) dx - \int_0^{t} \int_\Omega W_h \cdot \phi_t dx dt  \; 
$$
$$
=  -  \int_0^{t}  \nu\,[ W_h (\tau, \cdot) ,\phi (\tau, \cdot) ] dt \; + \;  \int_0^{t} \int_\Omega \bigg(- (u_* \cdot \nabla) W_h-(W_h \cdot \nabla) u_* \bigg) \phi dx d\tau$$
\begin{equation}\label{identityWh} + \;  \int_0^{t} \int_\Omega \bigg( F_h  \; +\, P_H (f_j(\tau,x;h)  - f_j (\tau, x;0)) \bigg) \phi dx d\tau, \;\; t \in [0, T^*]
\end{equation}
where  $F_h =  (w_h \cdot \nabla) (u_* - u_h)$. 
Then, making use of Remark \ref{ch} (see also calculations in Step 2 of subsection 4.2 below), we can derive, similar to (\ref{whest}) and \cite{Lad2}, page 167:
$$
\| W_h \|_{C ([0, T^*]; [\Lsp^2 (\Omega)]^2) \bigcap \Lsp^2([0, T^*; V)} 
$$
\begin{equation}\label{edelta} 
\leq C^*  \| f_j(\tau,x;h)  - f_j (\tau, x;0) \|_{[\Lsp_{2, 1}(Q_{T^*})]^2} +  C^*  \| F_h \|_{[\Lsp_{q_2, r_2 }(Q_{T^*})]^2}, 
\end{equation}
where $ q_2 = 2q/(q+1) = 4/3, r_2 = 2r/(r+1) = 4/3$ and $ C^*$ is from (\ref{est:wj}).
 \n
 In turn, see again (\ref{f*est}) in subsection 4.2 below:
$$
\| F_h \|_{[\Lsp_{q_2, r_2}(Q_{T^*})]^2}   \; \leq \; K \| \nabla u_* - \nabla u_h \|_{[[\Lsp^2 (Q_{T^*})]^2]^2} \|w_h \|_{[\Lsp^4 (Q_{T^*})]^2}
$$
\begin{equation}\label{eFh}
 \leq \; K_*  \| \nabla u_* - \nabla u_h \|_{[[\Lsp^2 (Q_{T^*})]^2]^2}  T^* \, {\rm meas}^{1/2} (\Omega) C_{\Omega, r},
 \end{equation}
where $ K_* > 0 $ is some constant and we used (\ref{hi2d}) and (\ref{4q}) to derive the 2nd inequality.

\n
{\bf Step 2.}
Note next  that, due to (\ref{eq:H2})-(\ref{elf}), Remarks \ref{fp} and \ref{nseest}, (see  (\ref{th11}) and (\ref{fje})), for some constant $k_r$, depending on $ r$, 
$$
\|  u_* -  u_h \|_{C([0, T^*]; H)} + \| \nabla u_* - \nabla u_h \|_{[[\Lsp^2 (Q_{T^*})]^2]^2}  \; \leq D^*  \| 0 \cdot f_{j}(\cdot,\cdot;0)- h f_{j}(\cdot,\cdot;h) \|_{[\Lsp_{2}  (Q_{T^*})]^2}
$$
$$
\leq \; D^*\sqrt{T^*} \, {\rm meas}^{1/2} \, (\Omega) \,   C_{\Omega, r} | h |
\rightarrow \; 0 \;\;\;\; {\rm as} \;\; h \rightarrow 0,
$$
$$
\| f_{j}(\cdot,\cdot;h)- f_{j}(\cdot,\cdot;0) \|_{[\Lsp_{2,1}  (Q_{T^*})]^2} \;  
\leq \; T^*  \, {\rm meas}^{1/2} \, (\Omega) \, \| f_{j}(\cdot,\cdot;h)- f_{j}(\cdot,\cdot;0) \|_{[\Lsp^\infty  (Q_{T^*})]^2}
$$
$$
\leq \; k_r  T^*  \, {\rm meas}^{1/2} \, (\Omega) \, \|z (\cdot;h)- z (\cdot;0) \|_{[C ([0, T^*]; \R^2)]^n}
\leq \; k_r (T^*)^{2}  \, {\rm meas} \, (\Omega) \, \|  u_* -  u_h \|_{C([0, T^*]; H)} 
$$
\begin{equation}\label{crate}
\leq \; D^* k_r (T^*)^{2.5}  \, {\rm meas}^{3/2} \, (\Omega) \,   C_{\Omega, r} | h | \; \rightarrow \; 0 \;\;\;\; {\rm as} \;\; h \rightarrow 0.
\end{equation}

\n
{\bf Step 3.} Estimates (\ref{edelta})- (\ref{crate})
  yield that
\begin{equation}\label{limwj}
\| w_h - {\partial u\over \partial v_j} |_{v_j's=0} \|_{C ([0, T^*]; H) \bigcap \Lsp^2 (0, T^*; V)} \;\leq \; {\cal C}(r, \Omega, T^*) | h | \;  \rightarrow 0 \;\;{\rm as} \; h \rightarrow 0,
\end{equation}
where ${\cal C}(r, \Omega, T^*) > 0 $ is  defined by $r, \Omega, T^*$ and  
\begin{equation}\label{c0}
{\cal C}(r, \Omega, T^*) \rightarrow \;\; {\rm as} \;  T^* \rightarrow 0.
\end{equation}
This   completes the proof of 
Lemma \ref{lem:wj}. $\diamond$

 \begin{rem}\label{ph}
 We would like to note here that the convergence rate in (\ref{limwj}) is linear with respect to $ | h |$.
 \end{rem}


\section{Derivatives  ${\partial z_i\over \partial v_j} |_ {v_j's = 0}$ as solutions to  Volterra equations: 2-$D$ case}\label{dzv}
We intend to show that 
$$
{\partial z_i\over \partial v_j} |_ {v_j's = 0} \; \mathop{=}^\Delta \;  \lim_{h\rightarrow 0}  \frac{z_i(t;h) - z_i(t; 0)}{h},
 $$
 where the limit is taken in $ [C[0,T^*]]^2$-norm exists.
To this end, we will use the integral form of equations  (\ref{eq:ode}):
$$
z_i(t;h) 
=z_{i,o}+\,{1\over \mes(S(0))}\,\int_0^t\int_{z_i(\tau; h)+S(0)} u(\tau,x;  h)\,dx\,d\tau
$$
$$
=z_{i,o}+\,{1\over \mes(S(0))}\,\int_0^t\int_{S(0)} u(\tau,x-z_i(\tau; h); h)\,dx\,d\tau.
$$
Respectively:
$$
 \frac{z_i(t;h) - z_i(t; 0)}{h} \; 
 $$
 $$
 = \,{1\over \mes(S(0))}\,\int_0^t\int_{S(0)} \frac{ u(\tau, x-z_i(\tau; h);h) - u(\tau,x-z_i(\tau; 0); 0)}{h} \,dx\,d\tau
 $$
$$
= \,{1\over \mes(S(0))}\,\int_0^t\int_{S(0)} \frac{u(\tau,x-z_i(\tau; h); h) - u(\tau,x-z_i(\tau; 0);h)}{h} \,dx\,d\tau \; 
$$
\begin{equation}\label{vj1}
+ \; {1\over \mes(S(0))}\,\int_0^t\int_{S(0)} \frac{u(\tau,x-z_i(\tau; 0); h ) - u(\tau,x-z_i(\tau; 0);0)}{h} \,dx\,d\tau. 
\end{equation}

\n
In the previous section we studied the integrand in the 2nd term (i.e., $ w_h$ and its limit properties as $ h \rightarrow 0$)), see Lemma \ref{lem:wj}.  

\bigskip
\n
{\bf Evaluation of the integrand in the 1st term on the right in (\ref{vj1}).} Due to assumption (\ref{conv}) and Remark \ref{fp}, without loss of generality (namely, for sufficiently small $h$), we can assume that for some $ \mu >0$ $z_i (t;h) \in B_\mu (z_i(t;0)) \subset \Omega,\; t \in [0, T^*]$  and 
$$
 \rho(s, t, x) \; \mathop{=}^ \Delta \; (1-s)(x-z_i (t;0)) + s (x - z_i (t;h)) 
 $$
 $$
 = x -  [(1-s)(z_i (t;0) + s z_i (t;h))] \in \Omega, s \in [0, 1], \; x \in S(0), \; t \in [0, T^*].
 $$
We claim that, due to assumption (\ref{eq:H2}) and Remarks  \ref{fp} and \ref{nseest}:
$$
u(t,x-z_i(t; h); h) - u(\tau,x-z_i(t; 0); h) 
$$\begin{equation}\label{vj2}
= \;  u_x ( t,x-z_i(t; 0); 0)   (z_i (t;0))  - z_i (t;h)) \; + \; G (t, x;h)  (z_i (t;0))  - z_i (t;h)),
\end{equation}
where  $u_x $ is the Jacobian matrix of the function $u (t, x)$ with respect to $x$ and 
\begin{equation}\label{vj5}
\| \int_{S(0) } G (t, x ;h) dx \|_{\R^2} \leq O(h) \| \nabla u(t, \cdot; h) \|_{[[\Lsp^2 (\Omega)]^2]^2}   \;\; {\rm as} \; h \rightarrow 0, \;\; t \in [0, T^*].
\end{equation}

\begin{rem}\label{O(t)}
Here and below, when we use a term like $ O(p)$ we assume that it may depends on the given parameters in the original problem (such as $ T^*, \Omega, r, u_0,  v_i$'s, selected indeces) but the limit property $O(p) \rightarrow 0 $ as $ p \rightarrow 0$ holds uniformly over such fixed parameters.
\end{rem}

\bigskip
\n
Indeed, for example, if $ u = (u_1, u_2), z_i = (z_{i1}, z_{i2})$, then:
$$
 u_1 (t, \rho (1, t, x);h) -  u_1 (t, \rho (0, t, x);h) = u_1 (t,x-z_i(t; h); h) - u_1 (t,x-z_i(t; 0); h)  \; 
$$
$$
= \;  u_{1x_1}  ( t,  \rho(s_1, t, x); h)   (z_{i1} (t;0))  - z_{i,1} (t;h)) \; + u_{1x_2}  ( t,\rho(s_1, t, x); h)   (z_{i2} (t;0))  - z_{i,2} (t;h)) ,
$$
where  $s_1 \in [0, 1]$ and point $ \rho(s_1, t, x)$ lies in the line interval connecting points $x - z_i(t; 0) $ and  $x- z_i(t; h) $ whose length tends to zero as $ h \rightarrow 0$ due to Remark  \ref{fp}. Then, 
$$
| \int_{S(0) } (u_{1x_1}  ( t,  \rho(s_1, t, x); h) - u_{1x_1}  ( t,  \rho(0, t, x); h)) \,dx | \
$$
$$
= \; | \int_{S(0) } (u_{1x_1}  ( t,  \rho(s_1, t, x); h) - u_{1x_1}  ( t,  x-z_i(t; 0); h)) \,dx | 
$$
$$
= \; | \int_{S(\rho(s_1, t, x) }  u_{1x_1}  ( t,  x ; h)    \,dx  \; 
-\;  \int_{S(z_i(t; 0) )}   u_{1x_1}  ( t,  x; h)  \,dx |
$$
$$
\leq \; | \int_{S(\rho(s_1, t, x) )\backslash S(z_i(t; 0) )} | u_{1x_1}  ( t,  x ; h) |  \,dx | \; 
+ \;   \int_{S(z_i(t; 0) )\backslash S(\rho(s_1, t, x))} |  u_{1x_1}  ( t,  x; h) | \,dx |
$$
$$
\leq \;  \mes (S(\rho(s_1, t, x))\backslash S(z_i(t; 0) ))^{1/2}  \|  u_{1x_1}  ( t,  \cdot  ; h) \|_{\Lsp^2 (\Omega)}  \; 
$$
$$
+ \;  \mes (S(z_i(t; 0) )\backslash S(\rho(s_1, t, x) )))^{1/2}  \|  u_{1x_1}  ( t,  \cdot  ; h) \|_{\Lsp^2 (\Omega)}. 
$$
In turn,  combining (\ref{eq:H2}) and Remark \ref{nseest}  yields  (\ref{vj5}).

\bigskip

\n
\underline{{\bf Volterra equations.}} Combining (\ref{vj1}) and (\ref{vj2})  yields the following Volterra equation for $  \psi_h = \frac{z_i(\cdot;h) - z_i(\cdot; 0)}{h} $:
\begin{equation}\label{veq}
((I  + {\mathbb A} + {\mathbb A_h}) \psi_h)(t) \; \mathop{=}^\Delta  \; \psi_h (t) +   \int_0^t \bigg( {\mathbb K}_0 (t, \tau) + {\mathbb K}_h (t, \tau) \bigg) \psi_h (\tau) d \tau =  g (t;h),
\end{equation}
$$
 I  + {\mathbb A} + {\mathbb A}_h: [C[0, T^*]]^2 \rightarrow [C[0, T^*]]^2, 
$$
where $ I $ is the identity operator and
\begin{equation}\label{vk}
{\mathbb K}_0 (t, \tau) \; = \; {-1\over \mes(S(0))}\, \int_{S(0)} u_x ( \tau,x-z_i(\tau; 0); 0)dx, \;\; {\mathbb K}_0 \in \Lsp^2 ((0, T^*) \times (0,T^*)),
\end{equation}
\begin{equation}\label{vk2}
{\mathbb K}_h (t, \tau) \; =  \;  {-1\over \mes(S(0))}\,\int_{S(0)}  G (\tau, x; h) dx, 
\end{equation}
 $$
g(t;h) \; =  \; {1\over \mes(S(0))}\,\int_0^t\int_{S(0)} {\partial u\over \partial v_j} |_{v_j's=0}  (\tau,x-z_i(\tau; 0 ))  \,dx\,d\tau \; +  \; H (t; h)$$
$$
H (t; h) \; = \; {1\over \mes(S(0))}\,\int_0^t\int_{S(0)} \frac{u(\tau,x-z_i(\tau; 0); h ) - u(\tau,x-z_i(\tau; 0);0)}{h} \,dx\,d\tau \;$$
\begin{equation}\label{Hest}
-  \; {1\over \mes(S(0))}\,\int_0^t\int_{S(0)} {\partial u\over \partial v_j} |_{v_j's=0}  (\tau,x-z_i(\tau; 0 ))  \,dx\,d\tau,  \; \| H(\cdot;h) \|_{C([0, T^*]; \R^2)} \rightarrow 0 \; {\rm as} \; h \rightarrow  0.
\end{equation}
The latter limit property is due to Lemma \ref{lem:wj}.
It is well-known that  operator $ I  + {\mathbb A} + {\mathbb A}_h$ in (\ref{veq}) is bijective and  has bounded inverse due to the Open Mapping Theorem. Recall that in (\ref{h1}) we assumed that $ | h | \leq 1$.

\bigskip
\n
{\bf Assumption on $ T^*$.}  {\em Recall that in (\ref{h1}) we assumed that $ | h | \leq 1$.
Without loss of generality, from now on, we can assume that $ T^*$ is small enough to ensure (making use of estimates in Remark \ref{nseest}) that}
\begin{equation}\label{anorm}
\| {\mathbb A} \| < \frac{1}{4}, \;\;\;\; \| {\mathbb A}_h \|  <  \frac{1}{4} \;\;\;\; \forall | h | \leq 1.
\end{equation}

\bigskip
\n
Respectively, in this case,  there exists a constant $ M_o $ such that
$$
\| \psi_h \|_{C([0, T^*]; \R^2)} \; \leq \; M_o, \;\; | h | \leq 1.
$$
Hence,  in view of  (\ref{vj5}), (\ref{Hest}) and (\ref{limwj}),
we can  pass to the limit in (\ref{veq}), described as
$$
\psi_h  =  \bigg(I  + {\mathbb A} \bigg)^{-1} \bigg(g(t;h) - {\mathbb A}_h \psi_h \bigg) \; = \; \bigg( \sum_{k = 0}^\infty {\mathbb A}^k \bigg) \bigg(g(t;h) - {\mathbb A}_h \psi_h \bigg),
$$
 in the $[C[0, T^*]]^2$-norm as $ h \rightarrow 0$ to obtain the existence of 
 $$ 
 {\partial z_i (t) \over \partial v_j} |_ {v_j's = 0} \; 
 =   \; \lim_{h\rightarrow 0}  \frac{z_i(t;h) - z_i(t; 0)}{h}
 $$ as the unique solution to the following limit Volterra equation: 
 $$
{\partial z_i (t) \over \partial v_j} |_ {v_j's = 0} \; 
 =  - \,\int_0^t \bigg\{ {1\over \mes(S(0))}\, \int_{S(0)}  u_x ( \tau,x-z_i(\tau; 0); 0)  dx  \bigg\}   \bigg[{\partial z_i (\tau) \over \partial v_j} |_ {v_j's = 0} \bigg] \,d\tau
 $$
\begin{equation}\label{vj8}
+ \; {1\over \mes(S(0))}\,\int_0^t\int_{S(0)} {\partial u\over \partial v_j} |_{v_j's=0}  (\tau,x-z_i(\tau; 0 ))  \,dx\,d\tau, 
\end{equation}
with
$$
\| {\partial z_i  \over \partial v_j} |_ {v_j's = 0} \|_{C([0, T^*]; \R^2)} \; 
$$
$$
\leq \;  \bigg( \sum_{k = 0}^\infty \| {\mathbb A}_t^k  \| \bigg) \; \| {1\over \mes(S(0))}\,
\int_0^{(\cdot)} \int_{S(0)} {\partial u\over \partial v_j} |_{v_j's=0}  (\tau,x-z_i(\tau; 0 ))  \,dx\,d\tau \|_{C([0, T^*]; \R^2)},
$$
\begin{equation}\label{zvnorm}
 \leq \; L_0 t^2    C_{\Omega, r},
 \end{equation}
where $ L_o > 0 $ is some constant,  $ {\mathbb A}_t$ is calculated as $ {\mathbb A} $ for the time interval $ (0, t) $ in place of $ (0, T^*)$ and we used (\ref{est:wj}).
Thus, we arrived at the following result.

\begin{lem}\label{eq:v} Assume (\ref{anorm}). Then derivatives 
$ {\partial z_i\over \partial v_j} |_ {v_j's = 0} , i = 1,  \ldots, n, j = 1, \ldots, 2n-3$ are elements of $C([0, T^*]; \R^2)$ and satisfy (\ref{vj8}).
\end{lem}

\n
Estimate (\ref{zvnorm})  immediately yields  the following lemma from (\ref{vj8}).

\begin{lem}\label{lem:estF1} Assume (\ref{anorm}). Then for any $ j = 1, \ldots, 2n-3$ and $t \in [0, T^*]$:
$$
\| {\partial z_i  \over \partial v_j} |_ {v_j's = 0}  \; - \;   {1\over \mes(S(0))}\,\int_0^{(\cdot)} \int_{S(z_i(\tau;0))}
{\partial u\over \partial v_j} |_{v_j's=0}  (\tau,x)  \,dx\,d\tau \|_{[C[0,t]]^2}
$$
\begin{equation}\label{est1}
=  \;  \;    t^2   O (t) \;\; {\rm as} \; t \rightarrow 0.
\end{equation}
\end{lem}


\section{Proofs of Theorems \ref{mm} and  \ref{thm:controllability}}\label{derf}

Let $w_j^o$ stand for solution to (\ref{eq:wj}) with  the right-hand side to be $ P_H f_j(0, \cdot; 0) $:
\begin{equation}\label{eq:wjo}
\left\{\begin{array}{ll}
w_{jt}^o \; = \; \nu\,\Delta w_j^o - (u_*\cdot\nabla) w_j^o \\
-(w_j^o \cdot\nabla) u_*+\, P_H f_j (\cdot, \cdot;0)  - \nabla p_j & \mbox{ in } (0,T^*)\times\Omega,\\
\Div \,w_j^o = 0& \mbox{ in } (0,T^*)\times\Omega, 
 \; {\rm i.e.}, \; w_j^o (t, \cdot)  \in H, \\
w_j^o=0 & \mbox{ in } (0,T^*)\times\partial\Omega,\\
w_j^o (0, \cdot)=0 & \mbox{ in } \Omega.
\end{array}\right.
\end{equation}
Then, due to (\ref{est:wj}), we have:
$$
\parallel {\partial u\over \partial v_j} |_{v_j's=0}  - w_j^o (t, \cdot) \parallel_{C([0, t]; H)}
$$
\begin{equation}\label{estw0}
 \leq \; 
C^* t \, {\rm meas}^{1/2} \, (\Omega) \,   \parallel  f_j(\cdot, \cdot; 0)  - f_j (0, \cdot;0) \parallel_{[L^\infty (Q_t)]^2} =  tO(t) \;\; {\rm as} \;  t \rightarrow 0,
\end{equation}
where we took into account (\ref{rot})-(\ref{elf}) and Remark \ref{nseest} in the last step.

\bigskip

\n
Combining (\ref{estw0}) with (\ref{est1}), yields:
$$
\| {\partial z_i  \over \partial v_j} |_ {v_j's = 0} \; -  \;{1\over \mes(S(0))}\,\int_0^{(\cdot)} \int_{S(z_i(\tau;0))} 
w_j^o \,dx\,d\tau \|_{[C[0,t]^2}  \;  
$$
\begin{equation}\label{estw02} 
= \; O(t) t^2  \;\; {\rm as} \;  t \rightarrow 0.
\end{equation}

\subsection{Proofs of Theorem \ref{mm} and  of Theorem \ref{thm:controllability} in the case of local controllability near  equilibrium}\label{thm:equil}

{\bf Step 1.} The equilibrium position for the swimmer in (\ref{eq:nse})-(\ref{eq:ode}) is the pair of solutions $ (u = 0 = u_0, z = z(0))$, initiated by 
the initial datum $ u_0 = 0 $, $ v_i = 0, i = 1, \ldots, 2n - 3$ and any set of $z_{i,0}, i = 1, \ldots, n$.
In this case (\ref{eq:wjo}) becames a system of {\em linear} nonstationary Stokes equations as follows:
\begin{equation}\label{eq:wjeq}
\left\{\begin{array}{ll}
w_{jt*}^o\; = \; \nu\,\Delta w_{j*}^o +\, P_H f_j (0, \cdot;0)  - \nabla p_j^0 & \mbox{ in } (0,T^*)\times\Omega,\\
\Div \,w_{j*}^o = 0& \mbox{ in } (0,T^*)\times\Omega, 
\; {\rm i.e.}, \; w_j^o  (t, \cdot)  \in H, \\
w_{j*}^o=0 & \mbox{ in } (0,T^*)\times\partial\Omega,\\
w_{j*}^o=0 & \mbox{ in } \Omega,
\end{array}\right.
\end{equation}
Respectively, its solution is represented by the following Fourier series \cite{Lad69}), \cite{KhBook}, Ch. 14:
$$
w_{j*}^o(t,x)  = \sum_{k = 1}^\infty \mathop{\int}_0^t e^{-\lambda_k (t - \tau)} \left(  \mathop{\int}_\Omega f_j^T (0, q;0)  \omega_k  dq  \right)  d \tau \;  \omega_k (x) 
$$
$$
=  t \sum_{k = 1}^\infty \frac{1- e^{-\lambda_k t}}{t\lambda_k} \left(  \mathop{\int}_\Omega (P_H f_j^T (0, \cdot;0) )(q) \omega_k dq   \right)
\omega_k (x),
$$
$$
=  \; t (P_H f_j^T (0, \cdot;0) )(x) \; - \; t \sum_{k = 1}^\infty \bigg( 1 - \frac{1- e^{-\lambda_k t}}{t\lambda_k} \bigg) \left(  \mathop{\int}_\Omega (P_H f_j^T (0, \cdot;0) )(q) \omega_k dq   \right)
\omega_k (x),
$$ 
where $^T$ stands for transposition. 

\n
{\bf Step 2.} Note now that   the function 
$
e(s) =  1 - \frac{1- e^{-s}}{s}, s > 0
$ 
tends to zero  as $ s\rightarrow 0+$ and to 1 as $ s \rightarrow \infty$ and is strictly monotone increasing on $ (0, \infty)$.  Therefore, 
\begin{equation}\label{PHF}
\| w_{j*}^o(t,\cdot) -   t P_H f_j^T (0, \cdot;0) \|_{[\Lsp (\Omega)]^2} \; = \;  t O(t) \;\;{\rm as} \; t \rightarrow 0,
\end{equation}
where $ O(t)$ depends on $ f_j (0, x;0)  $ (besides $ \lambda_k$'s, i.e., $ \Omega$), namely, on the rate of convergence of the series 
$$
\sum_{k = 1}^\infty  \left(  \mathop{\int}_\Omega (P_H f_j^T (0, \cdot;0) )(q) \omega_k dq   \right) \omega_k (x)
$$
in $ H$.
Combining (\ref{PHF}) and (\ref{estw02}) yields that the term $t P_H f_j^T (0, \cdot;0)$  will define the direction of vector ${\partial z_i  \over \partial v_j} |_ {v_j's = 0} $ as $t \rightarrow 0$, namely:
\begin{equation}\label{zijdir1}
\| {\partial z_i  \over \partial v_j} |_ {v_j's = 0} \;  -   \;{t^2 \over {2 \mes(S(0))}} \, \int_{S(z_i(\tau;0))} (P_H f_j^T (0, \cdot;0))(x) dx \|_{[C[0,t]]^2} \; = \; 
t^2 O(t) \;\; {\rm as}  \;  t \rightarrow 0+
\end{equation}
or
\begin{equation}\label{zijdir}
z_i (t) =  z_i(0) + {h t^2 \over {2 \mes(S(0))}} \, \int_{S(z_i(\tau;0))} (P_H f_j^T (0, \cdot;0))(x) dx + ht^2 \rho_j (t)  + h \zeta_j(h,t), \;t \in  [0, T^*], 
\end{equation}
where $  \| \rho_j \|_{[C[0,t]]^2}  = \;   O(t)$    and $  \| \zeta_j (h, \cdot) \|_{[C[0,t]]^2}  = \;  O(h)$.

\n
{\bf Step 3: Proof of Theorem \ref{mm}.} If we set in (\ref{h1})
\begin{equation}\label{h2}
v_j = h a_j \in \R, \;\; | h | \leq 1, \; \sum_{j = 1}^{2n - 3} a_i^2 = 1,
\end{equation}
we can repeat all the calculations leading to (\ref{zijdir}) with the force terms 
$$
 f(t,x) = f(t, x;h) = \sum_{j = 1}^{2n - 3}  v_j f_j (t,x) =  \sum_{j = 1}^{2n - 3} v_j f_j (t,x;  h), \;  
$$
in place  of the force term in (\ref{fh1}) and obtain (\ref{zijdirall}) instead. This proves Theorem \ref{mm}. $\diamond$

\n
{\bf Step 4.} Denote  
$$ 
V_{k,l,h} = \{v =  (v_k, v_l) \, | \, (v_k, v_l) = h(a_k, a_l), a^2 + a_l^2 = 1\}, \; V^{h_0}_{k,l} =  \bigcup_{0 \leq h \leq h_0} {\cal A}_{t,k,l} (V_{k,l,h}), \; h_0 \in [0, 1].
$$
Then, (\ref{zijdirall}) implies that the mapping $ {\cal A}_{t, k,l} : V^{h_0}_{k,l} \ni v \rightarrow z_i(t)$ can be represented as follows:
$$
z_i (t;h) \;  =  z_i(t;0) \; 
$$
$$
+   \; {t^2 \over { \mes(S(0))}} \, \bigg( v_k  \int_{S(z_i(\tau;0))} (P_H f_k^T (0, \cdot;0))(x) dx  \; + \;
v_l \int_{S(z_i(\tau;0))} (P_H f_l^T (0, \cdot;0))(x) dx  \bigg)  \; 
$$
\begin{equation}\label{zikldir}
+ \| v \|_{\R^2} t^2 \rho_{k,l}(t) + \| v \|_{|R^2} \zeta_{k,l} (\| v \|_{\R^2},t) , \;\; t \in  [0, T^*], 
\end{equation}
where $ \| \rho \|_{[C[0,t]]^2}  = O(t) $  and $   \| \zeta_{k,l} (\| v \|_{\R^2}, \cdot)  \|_{[C[0,t]]^2}  = O(\| v \|_{\R^2})$.

\n
Assuming that the vectors in (\ref{eq:indep_vect}) are linear independent and, due to Remark \ref{nseest}  (namely, on continuity of $ z_i$'s with respect to $ v_j f_j$'s), we can derive from  (\ref{zikldir}) that starting from some positive ``small'' $ h_0$ and for some $ T \in (0, T^*]$, the mapping 
$ {\cal A}_{T, k,l}$ is continuous, 1-1 and the range set $ {\cal A}_{T, k,l} (V^{h_0}_{k,l})$ is closed.
This also means that the images of the sets 
${\cal A}_{T, k,l} (V_{k,l,h}), h \in [0, h_0]$ will be closed curves encircling some neighborhoods of the point $ z_i(T;0) $ and   that  $ z_i(T;0) = z^* (T)$ is an internal point of the set 
$ {\cal A}_{T, k,l} (V^{h_0}_{k,l}) = \bigcup_{h \in [0, h_0]} {\cal A}_{T, k,l} (V_{k,l,h}) $, 
which implies the result of Theorem \ref{thm:controllability}  in the case of local controllability near equilibrium (as defined in \cite{KhBook}, Ch. 14). $
\diamond$


\subsection{Proof of Theorem \ref{thm:controllability}}\label{thm:gcase}

\n
{\bf Step 1.} We intend to prove Theorem \ref{thm:controllability} in the general case by adopting the formula (\ref{PHF}) to the {\em linear} system (\ref{eq:wjo}).
To this end,  we split solution to the latter into the sum of two functions:
$$
w_j^o = w_{j*}^o \; + \; u_e,
$$
where $u_e$ solves (\ref{eq:wjo}) in the case when $ P_H f_j (0, \cdot;0) = 0$, namely, for the following free term only (see also Remark \ref{rgal}):
 \begin{equation}\label{freeterm}
 f^* = (u_*  \cdot\nabla) w_j^*-(w_j^*\cdot\nabla) u_*,
\end{equation} 
\begin{equation}\label{eq:ue}
\left\{\begin{array}{ll}
u_{et}\; = \; \nu\,\Delta u_e +\, f^*  - \nabla p_j^* & \mbox{ in } (0,T^*)\times\Omega,\\
\Div \,u_e = 0& \mbox{ in } (0,T^*)\times\Omega, 
 \; {\rm i.e.}, \; u_e  (t, \cdot)  \in H, \\
u_e=0 & \mbox{ in } (0,T^*)\times\partial\Omega,\\
u_e=0 & \mbox{ in } \Omega,
\end{array}\right.
\end{equation}
Therefore, the general case of Theorem  \ref{thm:controllability} will follow as in subsection \ref {thm:equil}  if we will show that  
\begin{equation}\label{ueO}
\| u_e (t, \cdot) \|_{[\Lsp^2 (\Omega)]^2}  = t  O(t) \;\; {\rm as} \;\; t \rightarrow 0. 
\end{equation}

\bigskip
\n
{\bf Step 2.}  Once again, we invoke the results obtained within the proof of Theorems  1.1 in \cite{Lad2}, page 573 and Theorem 4.5 in \cite{Lad2}, page 166 establishing that 
(see ($4.8_1$)-($4.8_2$)  in \cite{Lad2}, page 156)
$$
f^* \in [\Lsp_{q_1, r_1} (Q_{T^*})]^2, \;\; q_1 = \frac{2q}{q+1} = \frac{4}{3},  \; r_1 = \frac{2r}{r+1}= \frac{4}{3},
$$
where $q=r=2$ are as in Remark \ref{ch}. This selection of $ q_1$  and $r_1$ satisfies the assumptions in this remark needed to apply the estimate (\ref{whest})/(\ref{est:wj}) to $ u_e$.
Thus, we obtain that:
\begin{equation}\label{ueest}
\| u_e \|_{C ([0, t]; [\Lsp^2 (\Omega)]^2) \bigcap \Lsp^2([0, t; V)}  \; \leq C_*  \| f^* \|_{[\Lsp^{4/3}(Q_{t})]^2}, \;\; t \in (0, T^*].
\end{equation}
In turn,  making use 
of H\"older's inequality, namely:  
\begin{equation}\label{hi2d}
| \int_{Q_t} \psi^{4/3}_1 \psi^{4/3}_2 dx dt | \; \leq \;  \bigg( \int_{Q_t} \psi^4_1  dx dt \bigg)^{1/3}  \bigg( \int_{Q_t} \psi^2_2  dx dt \bigg)^{2/3},  
\end{equation}
and, then, of the estimate  (\ref{4q}), applied to $ w_j^*$,  we can derive that for some positive constants $K, L$:
$$
\| f^* \|_{[\Lsp^{4/3} (Q_{t})]^2} \;
\leq \; K \bigg\{ \| u_* \|_{[\Lsp^4(Q_{t})]^2} \| \nabla w_j^* \|_{[[\Lsp^{2}(Q_{t})]^2]^2} 
+ \; \|  w_j^*  \|_{[\Lsp^4(Q_{t})]^2} \| \nabla  u_*  \|_{[[\Lsp^2(Q_{t})]^2]^2}  \bigg\} \; 
$$
\begin{equation}\label{f*est}
\leq \; K  L \bigg\{  \| u_* \|_{[\Lsp^4(Q_{t})]^2}  +  \| \nabla  u_*  \|_{[[\Lsp^2(Q_{t})]^2]^2}  \bigg\} \|  w_j^*  \|_{C ([0, T^*]; [\Lsp^2 (\Omega)]^2) \bigcap \Lsp^2([0, T^*; V)}  \; = \; t O(t) \;\; {\rm as} \;\; t \rightarrow 0.
\end{equation}
Here  
the last equality is due to estimates (\ref{whest2})/(\ref{est:wj}) applied to $w_j^*$.
Combining (\ref{ueest}) and (\ref{f*est}) yields (\ref{ueO}). This ends the proof of Theorem  \ref{thm:controllability}. $\diamond$


\section{{\bf 3-$D$ case. Proofs of Theorems \ref{mm} and  \ref{thm:c3d}.}} 
In the 3-$D$ case, we need to do the following adjustments in the above proofs relative to the 2-$D$ case.

\subsection{{\bf 3-$D$ case: Adjustments in Sections 2 and 3}}

\begin{itemize}
\item
Recall that that due to Theorem \ref{thm:wp23d},
\begin{equation}\label{u*est3d}
u_*, u_h, W_h  \in C ([0, T^*]; H) \bigcap C([0, T^*]; V)   \;\;  f_j (\cdot, \cdot; h) \in[ \Lsp^\infty (Q_{T^*})]^3.
\end{equation}
Due to  compact embedding $ V \subset  [H^1(\Omega)]^3 \subset [\Lsp^s (\Omega)]^3$  for $  s \in [1, 6) $(see, e.g.,  \cite{Lad69}, \cite{Lad2}, \cite{Brezis}).  Thus,
\begin{equation}\label{em3d1}
u_*, u_h, w_h  \in C([0, T^*]; V) \subset C([0, T^*]; [\Lsp_{s, \rho} (\Omega)]^3), \;  \rho > 0, s \in [1, 6).
\end{equation}
\item
Theorem 1.1 in \cite{Lad2}, pages 573 (see Remark  \ref{ch}) requires that the  squared 1-$D$ components of the 3-$D$ vector-function  $u_*$    and 1-$D$ components of the $3 \times 3$ matrix-function $ \nabla u_h$ (as the coefficients in (\ref{eq:wh})) are elements of the space $\Lsp_{q,r} (Q_{T^*})$, where 
$$
\frac{1}{r} + \frac{3}{2q} = 1, \;\;  q \in (1.5, \infty], \; r \in [1, \infty),
$$
while the free term $f_{j}(\cdot,\cdot;h) $ (we can ignore $- {1 \over h} \nabla (p (\cdot, \cdot; h) - p (\cdot, \cdot; 0))$ due to Remark \ref{rgal}) lies in $\Lsp_{q_1,r_1} (Q_{T^*})$, where \begin{equation}\label{q1r1}
 \frac{1}{r_1} + \frac{3}{2q_1} = 1 + \frac{3}{4}, \;  
q_1 \in [6/5, 2], \; r_1 \in [1, 2].
\end{equation} 

\n
We can select any suitable $ r_1, q_1$ in the above intervals, since $f_{j}(\cdot,\cdot;h) \in [\Lsp^\infty (Q_T^*))]^3$. Alternatively, we can select, e.g., $q_1 = 2, r_1 = 1$, to preserve the respective space for the free term in (\ref{whest}).

\n
For   the  squared 1-$D$ components of t  $u_*$ we can pick  $ q = 2 $ and $r =  4$, due  to (\ref{em3d1}).

\n
In view of (\ref{em3d1}), for the 1-$D$ components of  $ \nabla u_h$ we can also pick  $ q = 2 $ and $r = 4 $, emplying the embeding
$$
\nabla u_*, \nabla u_h \in C([0, T^*]; [[\Lsp^2 (\Omega)]^3]^3).
$$
\end{itemize}

\bigskip
\n
The above  implies that the results if the remainder of Section 3 and of Section 4  hold true in the 3-$D$ with the following corrections:
\begin{itemize}
\item
Constant  $C^*$ in (\ref{est:wj}) can be selected to be dependent only on $\| u_* \|_{C([0, T^*]; V)}$.
 \end{itemize}
 
\n

\subsection{{\bf 3-$D$ case: Adjustments in Section 4}}
The results of subsection 4.1 remain the same in the 3-$D$ case up to Step 4, which we can modify as follows.

\bigskip
\n
{\bf Section 4.1,  Step 4: 3-$D$-case.} 

\begin{itemize}
\item
Consider now a set of controls 
$$ 
V_{k,l,m, h} = \{ v = (v_k, v_l, v_m) = h (a_k, a_l, a_m)\, | \, a_k^2 + a_l^2 + a_m^2= 1\},  
$$
$$
V^{h_0}_{k,l} =  \bigcup_{0 \leq h \leq h_0} {\cal A}_{t,k,l,m} (V_{k,l,m, h}),  h_0 \in [0, 1].
$$
Then, assuming that the vectors 
\begin{equation}\label{eq:indep_vect3d}
 \int_{S(z_i(0))} (P_H  f_k (0, \cdot))(x)\,dx,  \; 
\int_{S(z_i(0))} (P_H  f_l (0, \cdot))(x)\,dx,  \; 
\int_{S(z_i(0))} (P_H  f_m (0, \cdot))(x)\,dx \end{equation}
are linear independent, as in Section 4.1 we can show  that for  some ``small'' positive $ h_0$ and $ T \in (0, T^*]$,  point $ z_i(T;0) = z^* (T)$ is an internal point of the set 
$\bigcup_{0 \leq h \leq h_0} {\cal A}_{T, k,l,m} (V_{k,l,m,h})$, which implies the result of Theorem \ref{thm:controllability}  in the case of local controllability near equilibrium in the 3-$D$ case. $
\diamond$
\end{itemize}

\n
{\bf Section 4.2: 3-$D$-case.} In subsection 4.2 we will need to make the following modifications in Step 2 to obtain (\ref{ueO}):\begin{itemize}
\item
Once again, we invoke the results obtained within the proof of Theorems  1.1 in \cite{Lad2}, page 573 and Theorem 4.5 in \cite{Lad2}, page 166 establishing that 
(see ($4.8_1$)-($4.8_2$)  in \cite{Lad2}, page 156)
$$
f^* =  (u_*  \cdot\nabla) w_j^* - (w_j^*\cdot\nabla) u_*  \in [\Lsp_{q^*_1, r^*_1} (Q_{T^*})]^3, \;\; q^*_1 = \frac{2q}{q+1} = \frac{4}{3},  \; r^*_1 = \frac{2r}{r+1} = \frac{8}{5} 
$$
and condition (\ref{q1r1}) holds for these $ q^*_1$ and $r^*_1$ (in place of  $q_1, r_1$) with the above-selected $q = 2$ and $r = 4$.

\n
Next, due to Lemma 1.1 in  \cite{Lad2}, pages 59-60  (on the space dual of  $\Lsp_{q_1, r_1}  (Q_t)$) and estimates (1.11)-(1.12) on page 137 in \cite{Lad2}, we have:
$$
\| f^* \|_{[\Lsp_{q^*_1, r^*_1} (Q_t)]^3} 
$$
\begin{equation}\label{f*e3d}
\leq \; M_o \bigg\{ \| u_*^2 \|^{1/2}_{[\Lsp_{q, r} (Q_t)]^3} \;\| \nabla w_j^* \|_{[[\Lsp^2(Q_t)]^3]^3} \; 
\; + \| w_j^* \|_{[\Lsp_{\bar{q}, \bar{r}} (Q_t)]^3} \;\| \nabla u_* \|_{[[\Lsp_{q,r}(Q_t)]^3]^3} \bigg\}
\end{equation}
for some $ M_o > 0$, where
$$
q = \frac{\bar{q}}{\bar{q} - 2} = 2, \bar{q} = 4 \;\; {\rm and} \;\; r = \frac{\bar{r}}{\bar{r} - 2} = 4, \bar{r} = \frac{8}{3}.
$$
\end{itemize}

\begin{rem}\label{l43d}
Let us recall  estimate  (3.4) in  \cite{Lad2}, page 75 for 3-$D$ case, namely:
\begin{equation}\label{4q3d}
\| \psi \|_{\Lsp_{q_*,r_*}(Q_{t})} \; \leq \; \beta \bigg( \max_{\tau \in [0, t]} \| \psi (\tau, \cdot) \|_{\Lsp^2(\Omega )}  + \| \nabla \psi \|_{[\Lsp^2(Q_{t})]^2} \bigg),
\end{equation}
$ 1/r_* + 3/(2q_*) = 3/4, r_* \in [2, \infty), q_* \in [2, 6]$.
\end{rem}

\begin{itemize}
\item
Due to (\ref{4q3d}), 
$$
\| w_j^* \|_{[\Lsp_{\bar{q}, \bar{r}} (Q_t)]^3} \; = \; \| w_j^* \|_{[\Lsp_{4, 8/3} (Q_t)]^3} \; \leq \; \beta_* \| w_j^*\|_{C ([0, T^*]; [\Lsp^2 (\Omega)]^2) \bigcap \Lsp^2([0, T^*; V)}
$$
for some $ \beta_* > 0$. Hence, for some $ M_o^* > 0$:
$$
\| f^* \|_{[\Lsp_{q^*_1, r^*_1} (Q_t)]^3} 
$$
\begin{equation}\label{f*e3d2}
\leq \; M^*_o \bigg\{ \| u_*^2 \|^{1/2}_{[\Lsp_{q, r} (Q_t)]^3} \;\| \nabla w_j^* \|_{[[\Lsp^2(Q_t)]^3]^3} \; 
\; +  \| w_j^* \|_{\Lsp^2  (0, T^*; V)} \;\| \nabla u_* \|_{[[\Lsp_{q,r}(Q_t)]^3]^3} \bigg\}.
\end{equation}
\item
Estimate (\ref{f*e3d2}) and the 3-$D$ version of estimate (\ref{est:wj}), applied for $ w_j^*$, yields that, instead of (\ref{ueest}), we have:
\begin{equation}\label{ueest3d}
\| u_e \|_{C ([0, t]; [\Lsp^2 (\Omega)]^3) \bigcap \Lsp^2([0, t; V)} \leq \hat{C}  \bigg\{ \| f_1^*  \|_{[\Lsp_{q_1^*, r_1^*}(Q_{t})]^3} \bigg\}= t O(t) \; {\rm as} \; t \rightarrow 0.
\end{equation}
This ends the proof of Theorem  \ref{thm:c3d}. 
\end{itemize}
\hfill $\diamond$


\section{Illustrating examples}
Let us recall some results from \cite{Kh7} and \cite{KhBook} (Ch. 13) allowing to calculate the averaged projections $\int_{S(z_i(0))} (P_H  f_k (0, \cdot, z(0)))(x)\,dx, k = 1, \ldots 2n-3$ in the cases when $ S(0) $ is a rectangle or a disc.

\begin{figure}[h]\label{fig:swim2prd}
\centering
\includegraphics[scale=0.6]{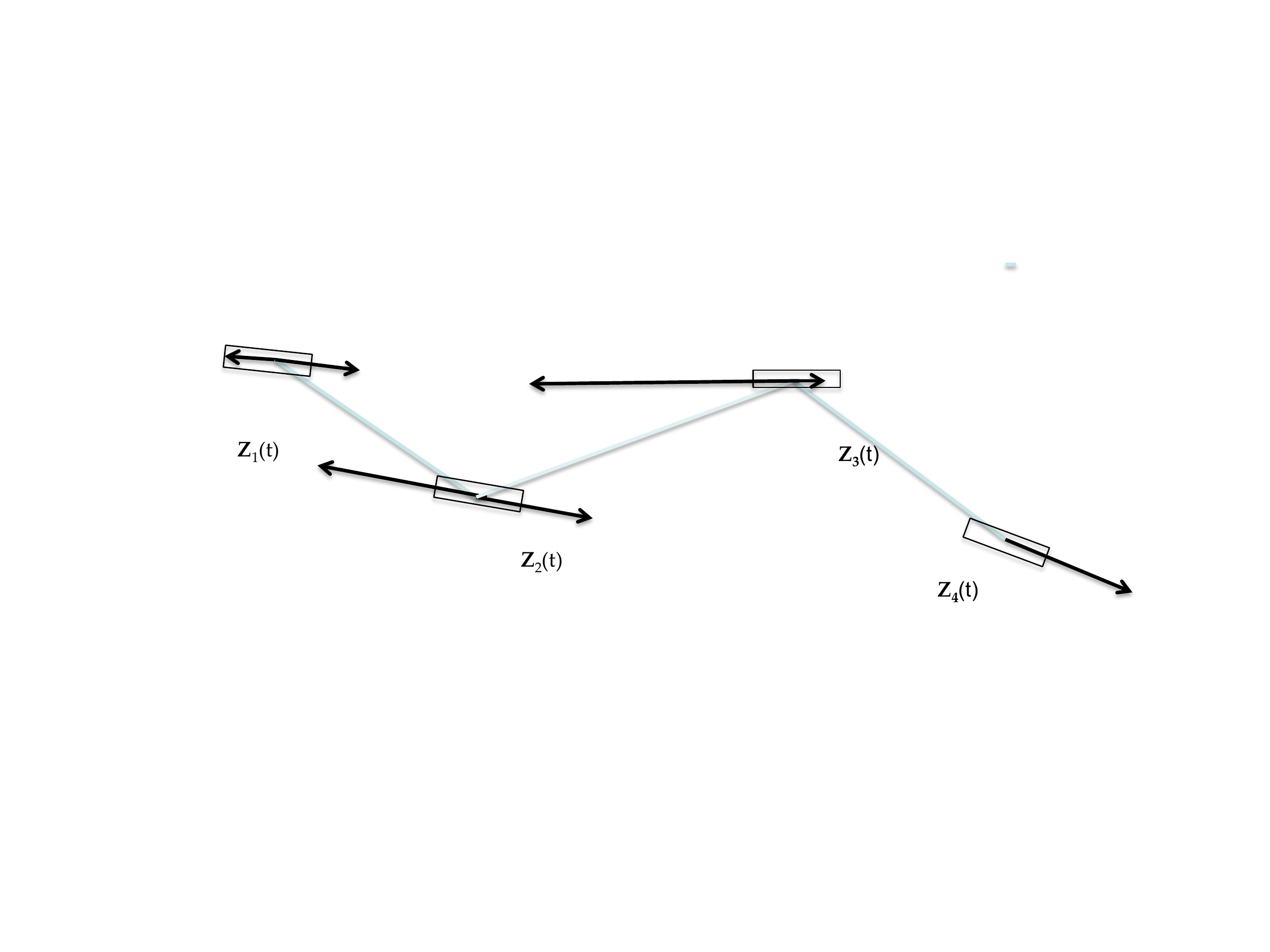} 
\caption{The (approximate) forces acting on the 2-D swimmer from Figure 2  in  the fluid.}
\end{figure}
\n
Assume that
\begin{equation}\label{eq:rt}
S (0) = S_0  = \{(x = (x_1, x_2) \mid -p < x_1 < p, \; -q < x_2 < q \},
\end{equation}
where $ p $ and $q $ are ``small'' positive numbers.

\bigskip
\begin{thm}[cited from \cite{Kh7} and \cite{KhBook} (Ch. 13)]\label{pt}  Let $ b = (b_1, b_2) $ be a given 2-$D$ vector and $ S (0) = S_0$ as in (\ref{eq:rt}) lie in $ \Omega $. Let $ q, p, q^{1-a}/p \rightarrow 0+$ for some  $ a \in (0, 1)$. Then
\begin{equation}\label{eq:pg}
\frac{1}{  {\rm mes} \, \{ S (0) \}}
\mathop{\int}_{S (0)} (P b \xi)(x) dx \; = \; (b_1,  0)  \;  + \; O(q^a ) \; + \; O (q^{1-a}/p) \; + \; O (p)
\end{equation}
as $ q, p, q^{1-a}/p \rightarrow 0+$, where $ \xi (x) $ is the characteristic function of $ S(0)$.
\end{thm}

\n
We can interpret this theorem as that the average projection of a force  $ b = (b_1, b_2) $, acting on a small narrow rectangle, in the fluid velocity space is approximately equal to its projection on the direction parallel to the longer side of the rectangle.

\bigskip
\n 
Figure 3 shows the transformation of internal forces of the swimmer from Figure 2 into the forces which actually interact with surrounding medium when the swimmer is in the fluid.

\begin{thm}[cited from \cite{Kh7} and \cite{KhBook} (Ch. 13)]\label{dt} Let $ b = (b_1, b_2) $ be a given 2-$D$ vector and $ S (0) = S_0$ be a \underline{disc} of radius $r$ lying  in $ \Omega $. Then 
\begin{equation}\label{eq:dg}
\frac{1}{  {\rm mes} \, \{ S (0) \}}
\mathop{\int}_{S (0)} (P b \xi)(x) dx \; = \; \frac{1}{2} (b_1,  b_2)  \;  + \; O(r ) 
\end{equation}
as $ \; r \rightarrow 0+$, where $ \xi (x) $ is the characteristic function of $ S(0)$.
\end{thm}

\n
We can interpret this theorem as that the average projection of a force  $ b = (b_1, b_2) $, acting on a small disc in the fluid velocity space, is approximately equal to its half.

\bigskip
\begin{lem}[cited from  \cite{Kh7} and  \cite{KhBook} (Ch. 13)]\label{rf} Let $ b = (b_1, b_2) $ be a given 2-$D$ vector.  Let $ S (0) \subset \Omega $ be a nonzero measure set  which  is strictly separated from $ \partial \Omega $  and   lies in an $ r$-neighborhood ($r>0$) of the origin. Then for any subset $ Q $ of $ \Omega$ of positive measure of diameter $ 2r$ (that is, it fits some ball of radius $ r$) which lies outside of some, say, $ d$-neighborhood ($ d > 0$) of $ S (0)$ and is strictly separated from $ \partial \Omega$ we have:
$$
\frac{1}{  {\rm mes} \, \{ Q \}}
\mathop{\int}_{Q} (P b \xi)(x) dx \; = \;  O(r)
\eqno(3.1)$$
as $ r \rightarrow 0+$, where $ \xi (x) $ is the characteristic function of $ S(0)$.
\end{lem}
\begin{rem}[Influence of ``remote body forces'']\label{rrf} 
We can interpret Lemma \ref{rf} as that  the effect of the force  $ b \xi (x) $ on  similar sized sets outside of its support $ S (0)$ is ``small'' if the size of $ S(0)$ is ``small''. In other words, the results of actions of swimmer's internal forces applied not directly to the body part at hand are ``negligible'' relative to the result of the forces applied directly on this body part.
\end{rem}

\bigskip
\n
We will use the above-cited results to illustrate possible applications of Theorems \ref{thm:controllability} and \ref{thm:ccm} as presented on  Figures 3-7.

\bigskip
\n
{\bf Examples 6.1: Local controllability of the center of mass/self-propulsion.}  Consider the  swimmer from  Figures 1-3  and assume that the rectangles forming its body are asymptotically small and satisfy the assumptions of  Theorem  \ref{pt},  and are oriented as on Figure 4-5. Then it follows from Theorem \ref{thm:ccm} that the swimmer on Figures 4-5 is locally controllable near  its center of mass by varying  various  pairs of $ (v_l, v_k)$ as long as condition  (\ref{eq:indep_vectc}) holds.

\n
For example, one can activate only the  pair of controls $ (v_3, v_5)$ in (\ref{elf2}) defining the elastic  forces acting between $z_1(t) $ and  $z_2(t)$ and between $z_3(t) $ and  $z_4(t)$, see Figures 4 and 5. In this case  the 1st pair of internal forces will result in an averaged  projected force acting on rectangle  $ z_1(t)$ approximately parallel to its  longer side (the elastic force acting on $z_2(t)$ can be ``neglected'' as perpendicular to the longer side  of this rectangle). In turn,  the 2nd pair of internal forces  will create a pair of averaged projected force acting on $ z_3(t)$ and $z_4(t)$ defined by  respective spatial orientations of these rectangles. The sum of these forces  is not co-linear (under our assumptions) to the  averaged  projected force acting on $ z_1(t)$,  see Figure 5. Thus, condition (\ref{eq:indep_vectc}) holds. In this example the swimmer is capable of local self-propulsion  (locomotion) -  moving of its center of mass under the actions of its internal forces.

\begin{figure}[h]\label{fig:swim2prd3}
\centering
\includegraphics[scale=0.6]{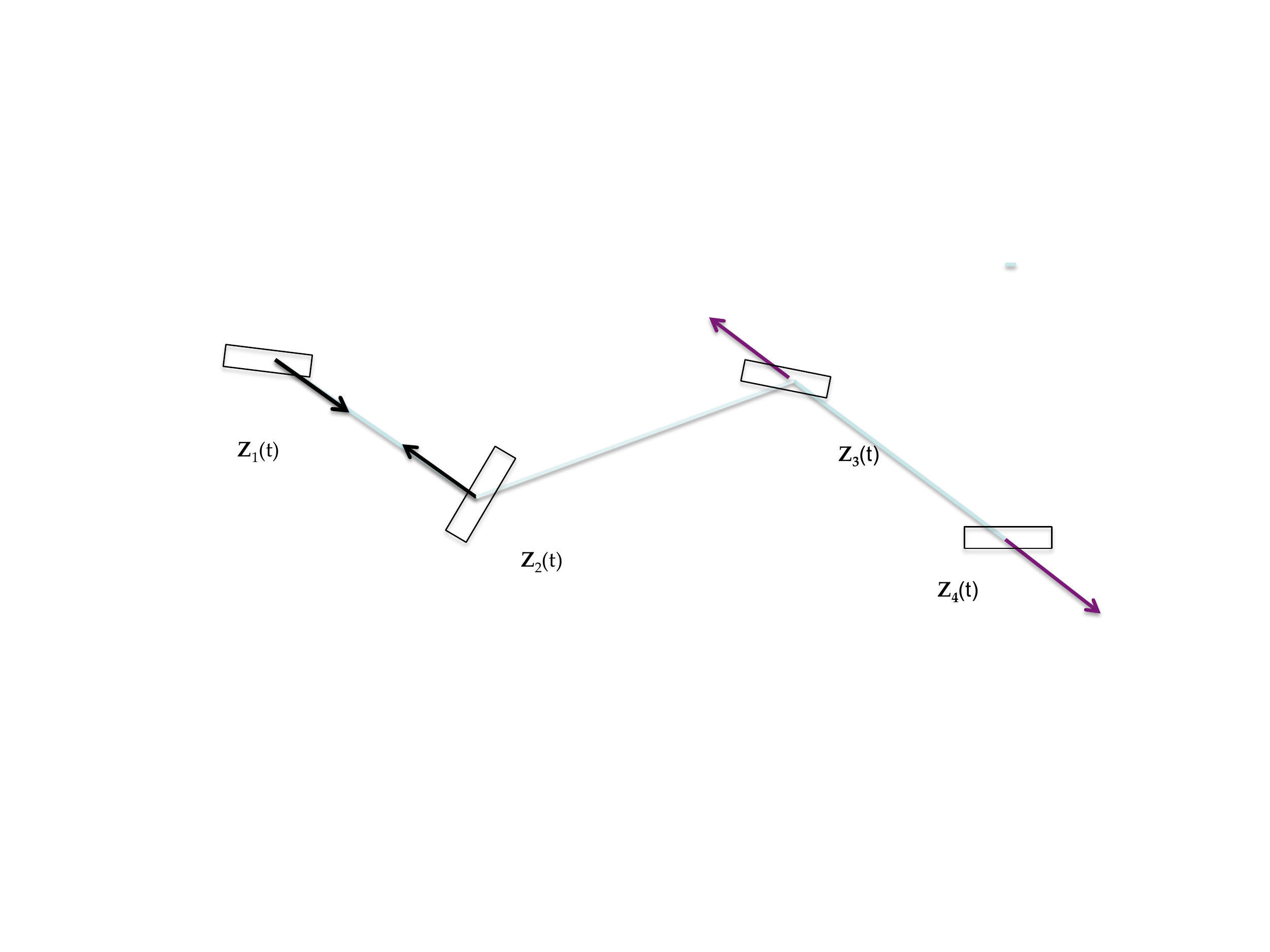} 
\caption{The 2-D swimmer from Figure 2  with 2 pairs of elastic forces acting.}
\end{figure}

\begin{figure}[h]\label{fig:swim2prd4}
\centering
\includegraphics[scale=0.6]{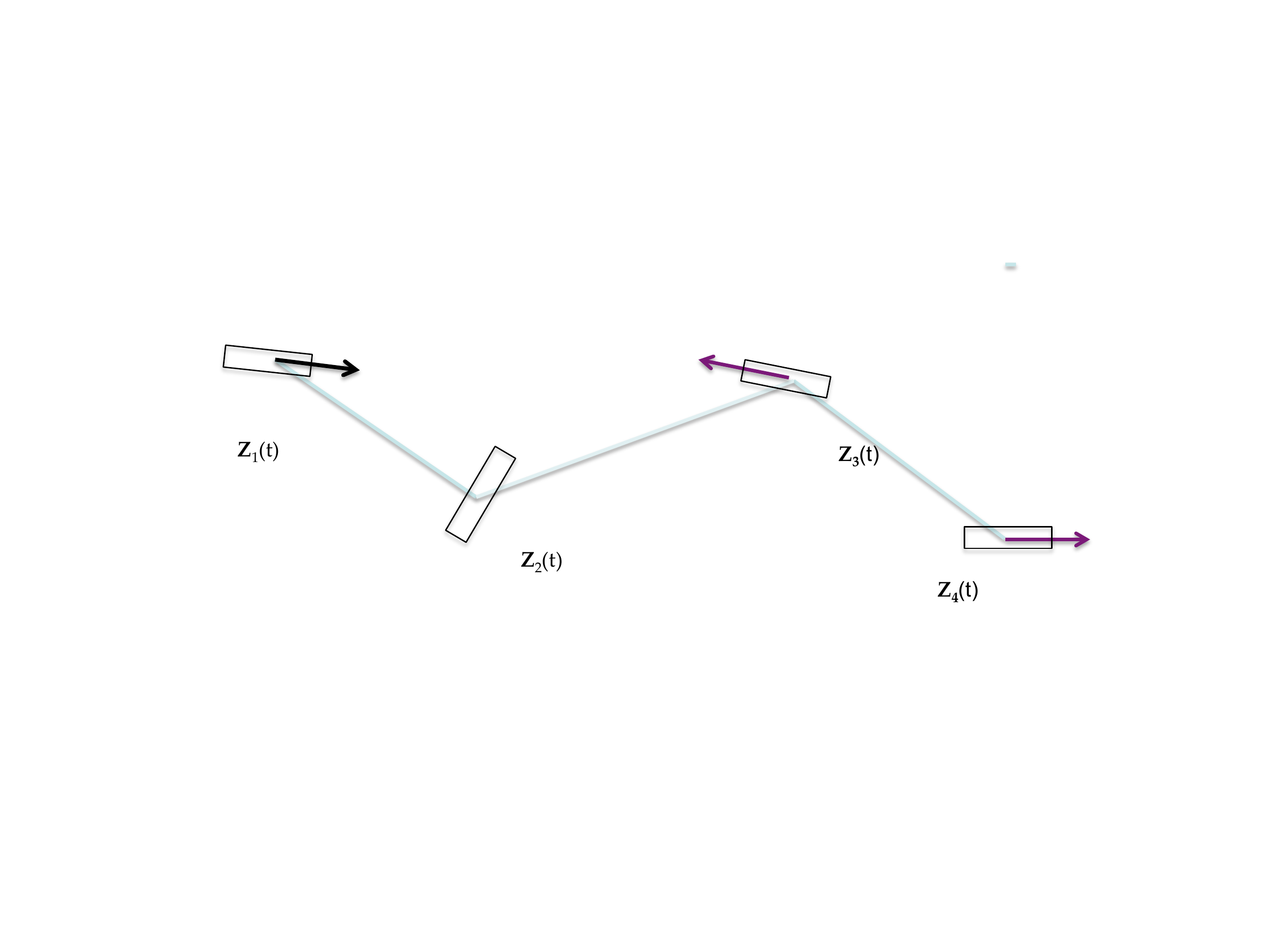} 
\caption{The  2-D swimmer from Figure 4  in  the fluid. Approximate averaged projected forces on the fluid velocity space are shown.}
\end{figure}

\begin{figure}[h]\label{fig:swim4d}
\centering
\includegraphics[scale=0.6]{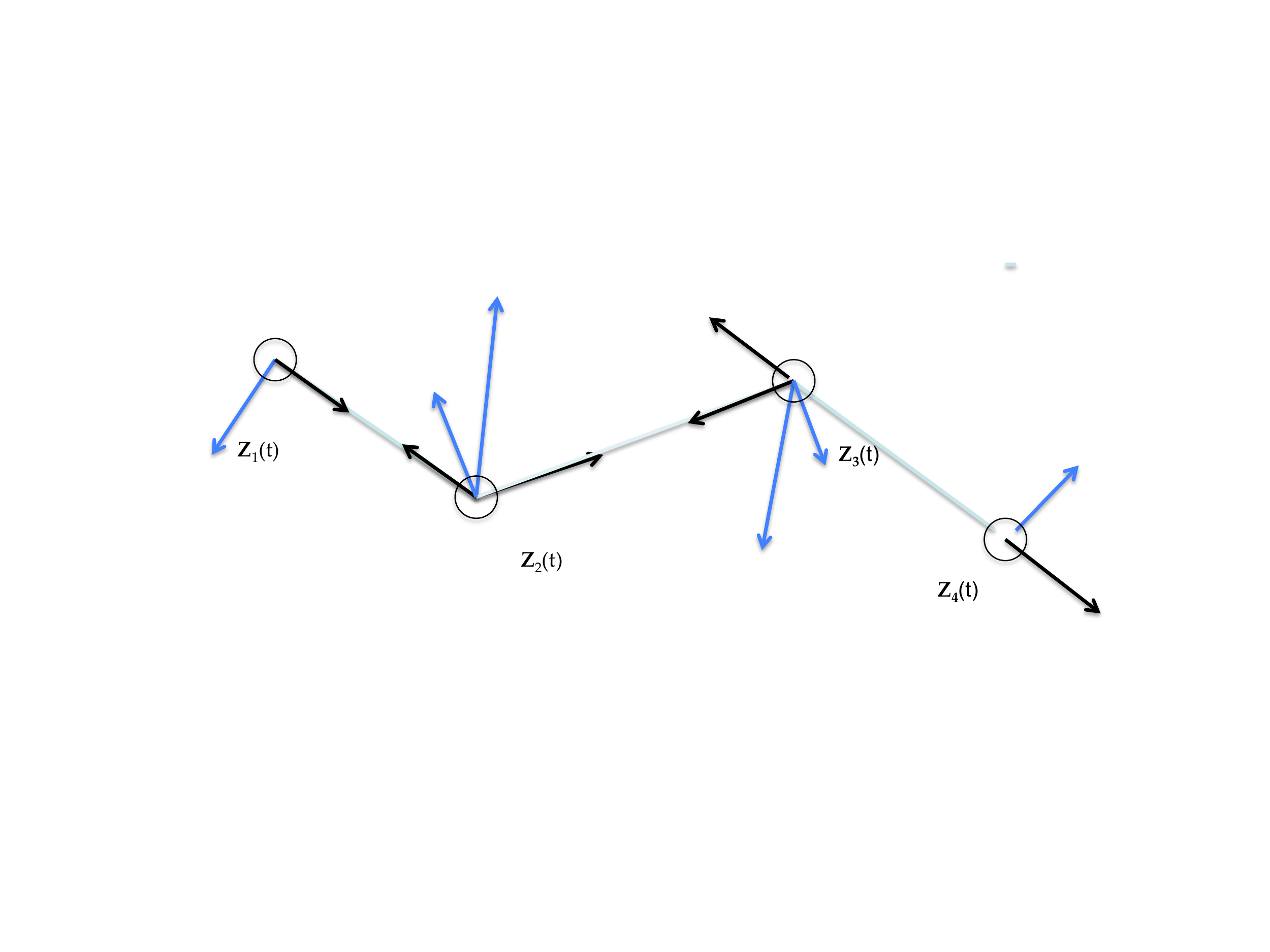} 
\caption{2-D swimmer consisting of 4 discs.}
\end{figure}

\begin{figure}[h]\label{fig:swim4d}
\centering
\includegraphics[scale=0.6]{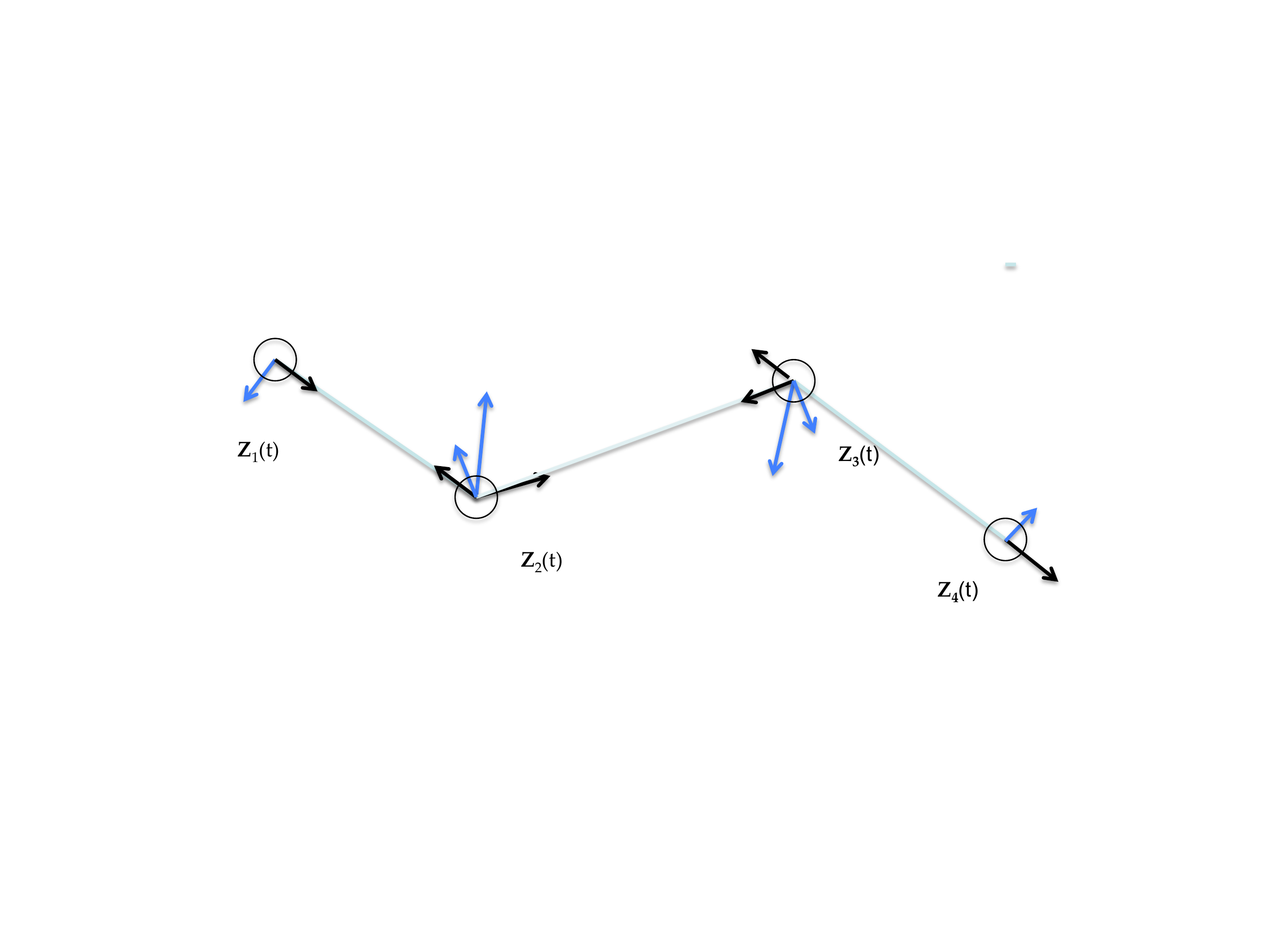} 
\caption{2-D swimmer from Figure 6 in the fluid. Approximate  averaged projection  forces on the fluid velocity space are shown.}
\end{figure}

\begin{figure}[h]\label{fig:swim3d}
\centering
\includegraphics[scale=0.6]{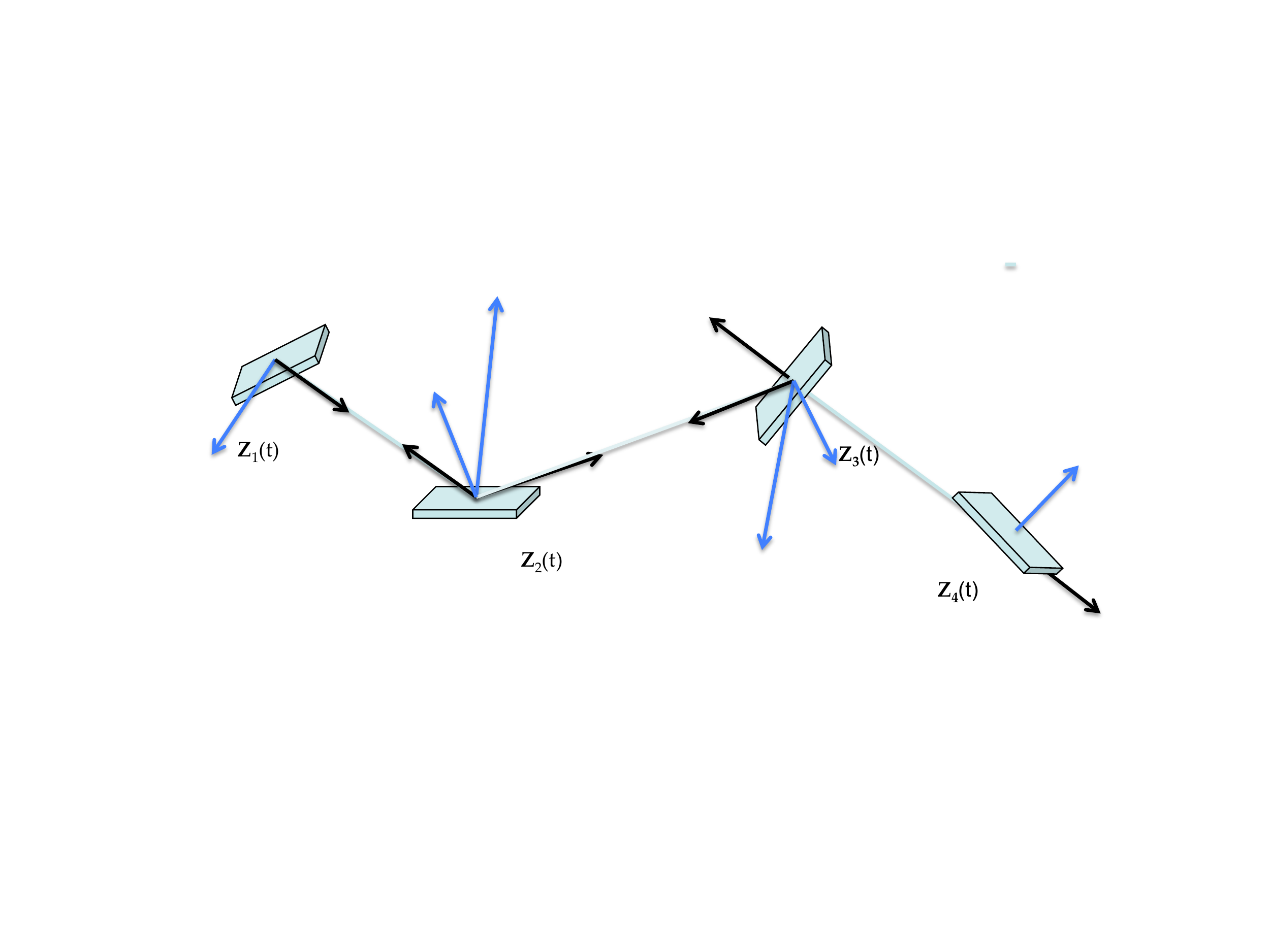} 
\caption{3-D swimmer consisting of 4 parallelepipeds.}
\end{figure}

\n
{\bf Remark on global controllability in Example 6.1.} The local controllability of the center of mass of the swimmer in Example 6.1 means that it can  move its center of mass within {\em some neighborhood}  in any direction by varying controls $ (v_3, v_5)$ in (\ref{elf2}). It does not mean  that all points $ z_i(t)$ will move in the same direction. For example,  in the case of shown on Figure 5 $ z_1 $ and $ z_3$ will (approximately)  move to the right, $ z_3$ will  move to the left, while $ z_2$ will not move.  A way to achieve the global swimming controllability of swimmer in Example 6.2 - {\em as a principal possibility to move the center of swimmer's mass between  any two points in $ \Omega$} - can be a combination of subsequent employment of various pairs of controls $ v_i$s which would move the swimmer in small increments towards the target point, while preserving its prescribed structure (such as maintaining  allowed limits for deviations of distances between points $ z_i$'s). This method was applied in \cite{KhBook}, Ch. 15 for a 2-$d$ swimmer whose body consisted of three rectangles and for the fluid described by the nonstationary Stokes equations (see also \cite{Kh4} for the 3-$D$ case). 

\bigskip
\n
{\bf Examples 6.2: On lack of self-propulsion in the case when the swimmer is formed by a set of discs.}  On Figures 6 and 7 we have the same configuration of a swimmer with the same internal forces but composed of small identical discs. Theorem \ref{dt}  implies that  in this case Theorem \ref{thm:ccm}  {\em does not guarantee  the self-propulsion of the  swimmer},  because the sum of all its averaged  projected forces on the fluid velocity space (responsible for the motion of the center of its mass)   will remain zero at all times.  Namely, these forces will preserve the directions of the original internal forces, while their magnitudes will be reduced by the same factor. In particular, both vectors in condition (\ref{eq:indep_vectc}) become zero-vectors. However,  we can  have the local swimmer's controllability near all  points $ z_i, i = 1, 2, 3, 4$, due to Theorem  \ref{thm:controllability}.

\n
The center of swimmer's mass can, in general, also  change its location in this example under certain circumstances, for example:
\begin{itemize}
\item
due to fluid's drifting motion (such as its natural ``flow'' associated with given  $ u_0$) or 
\item
due to fluid's turbulence, induced by the movements of swimmer's body parts inflicted by actions of its internal forces. 
\end{itemize}
\bigskip
\n
{\bf 3-$D$ examples.} Based on the results of \cite{Kh3}, extending Theorems \ref{pt} and \ref {dt} to the 3-$D$ case, Examples 6.1-2 can be modified  as follows:
\begin{itemize}
\item
Relation analogous to (\ref{eq:dg}) holds in  3-$D$  incompressible fluids with factor  1/3 when the discs are replaced by asymptotically small spheres or cubes, see  \cite{Kh3}.  Respectively, Example 6.2 will remain unchanged in the case when the body of a  3-$D$ swimmer consists of  any finite number of identical spheres or cubes.
\item
The results as in Example 6.1 can be obtained in the 3-$D$ case for a swimmer whose body consists of parallelepipeds (see Figure 8) whose proportions satisfy certain asymptotic assumptions qualitatively similar to those in Theorem \ref{pt}, see \cite{Kh3} and illustrating examples in \cite{Kh4}.
\end{itemize}



\n


\end{document}